\documentclass[11pt]{article}
\usepackage{mathrsfs}
\usepackage{amsthm}
\usepackage{amssymb}
\usepackage{amsmath}
\usepackage{graphicx}
\usepackage{color}
\usepackage{amsfonts}
\usepackage{float}
\usepackage{cite}
\usepackage[text={140mm,210mm},left=35mm,vmarginratio=1:1]{geometry}
\newtheorem{theorem}{Theorem}[section]

\newtheorem{lemma}[theorem]{Lemma}

\newtheorem{corollary}[theorem]{Corollary}

\numberwithin{equation}{section}
\normalsize

\begin{document}
\title{\textbf{The survival probability of the high-dimensional contact process with random vertex weights on the oriented lattice}}

\author{Xiaofeng Xue \thanks{\textbf{E-mail}: xfxue@bjtu.edu.cn \textbf{Address}: School of Science, Beijing Jiaotong University, Beijing 100044, China.}\\ Beijing Jiaotong University}

\date{}
\maketitle

\noindent {\bf Abstract:} This paper is a further study of Reference
\cite{Xue2015}. We are concerned with the contact process with
random vertex weights on the oriented lattice. Our main result gives
the asymptotic behavior of the survival probability of the process
conditioned on only one vertex is infected at $t=0$ as the dimension
grows to infinity. A SIR model and a branching process with random
vertex weights are the main auxiliary tools for the proof of the
main result.

\quad


\section{Introduction}\label{section one}

In this paper we are concerned with the contact process with random
vertex weights on the oriented lattice $\mathbb{Z}_+^d$ for $d$
sufficiently large, where $\mathbb{Z}_+=\{0,1,2,\ldots\}$. This
paper is a further study of Reference \cite{Xue2015}, which deals
with the critical value of the aforesaid process. First we introduce
some notations and definitions. For $x=(x_1,\ldots,x_d)\in
\mathbb{Z}_+^d$, we define
\[
\|x\|=\sum_{j=1}^dx_j
\]
as the $l_1$ norm of $x$. For $1\leq j\leq d$, we use $e_j$ to
denote the $j$th elementary unit vector of $\mathbb{Z}_+^d$, i.e.,
\[
e_j=(0,\ldots,0,\mathop 1\limits_{j \text{th}},0,\ldots,0).
\]
We use $O$ to denote the origin of $\mathbb{Z}_+^d$. For $x,y\in
\mathbb{Z}_+^d$, we write $x\rightarrow y$ when and only when
\[
y-x=e_j
\]
for some $j\in \{1,2,\ldots,d\}$.

For later use, we define a order $\prec$ on $\mathbb{Z}^d_+$. For
$x=(x_1,\ldots,x_d),y=(y_1,\ldots,y_d)\in \mathbb{Z}^+_d$, $x\prec
y$ when and only when there exists $j\in \{1,2,\ldots,d\}$ such that
$x_j<y_j$ while $x_i=y_i$ for any $i<j$.

Let $\rho$ be a random variable such that $P\big(\rho\in
[0,M]\big)=1$ for some $M\in (0,+\infty)$ and $P(\rho>0)>0$, then we
assign an independent copy $\rho(x)$ of $\rho$ on each vertex $x\in
\mathbb{Z}_+^d$. $\rho(x)$ is called the vertex weight of $x$. We
assume all these vertex weights are independent. After the vertex
weights are given, the contact process $\{C_t\}_{t\geq 0}$ on
$\mathbb{Z}_+^d$ with vertex weights $\{\rho(x)\}_{x\in
\mathbb{Z}_+^d}$ is a continuous time Markov process with state
space
\[
X=\big\{A:~A\subseteq \mathbb{Z}_+^d\big\}
\]
and transition rates function given by
\begin{equation}\label{equ 1.1 transition rate}
C_t\rightarrow
\begin{cases}
C_t\setminus\{x\} & \text{~at rate~} 1 \text{~if~}x\in C_t,\\
C_t\bigcup\{x\} & \text{~at
rate~}\frac{\lambda}{d}\sum\limits_{y:~y\rightarrow
x}\rho(x)\rho(y)1_{\{y\in C_t\}} \text{~if~}x\not \in C_t,
\end{cases}
\end{equation}
where $\lambda$ is a positive constant called the infection rate
while $1_A$ is the indicator function of the event $A$.

Intuitively, the process describes the spread of an epidemic on
$\mathbb{Z}_+^d$. Vertices in $C_t$ are infected while vertices out
of $C_t$ are healthy. An infected vertex waits for an exponential
time with rate one to become healthy while a healthy vertex $x$ may
be infected by an infected vertex $y$ when and only when
$y\rightarrow x$. The infection occurs at rate proportional to the
product of the weights on these two vertices.

The (classic) contact process is introduced by Harris in
\cite{Har1974}, where $\rho\equiv1$ and infection occurs between
nearest (un-oriented) neighbors. For a detailed survey of the
classic contact process, see Chapter 6 of \cite{Lig1985} and Part 1
of \cite{Lig1999}.

The contact process with random vertex weights is first introduced
in \cite{Pet2011} on the complete graph $K_n$ by Peterson, where a
phase transition consistent with the mean-field analysis is shown.
In detail, infected vertices dies out in $O(\log n)$ units of time
with high probability when $\lambda<\frac{1}{E(\rho^2)}$ or survives
for $\exp\{O(n)\}$ units of time with high probability when
$\lambda>\frac{1}{E(\rho^2)}$. In \cite{Xue2015}, Xue studies this
process on the oriented lattice and gives the asymptotic behavior of
the critical value of the process as the dimension $d$ grows to
infinity. When $P(\rho=1)=p=1-P(\rho=0)$ for some $p\in (0,1)$, the
model reduces to the contact process on clusters of the site
percolation, which is a special case of the model introduced in
\cite{Ber2011} with $n=1$. In \cite{Ber2011}, Bertacchi, Lanchier
and Zucca study the contact process on $G\times K_n$, where $G$ is
the infinite open cluster of the site percolation while $K_n$ is the
complete graph with $n$ vertices. Criteria judging whether the
process survives is given.

If the i.i.d. weights are assigned on the edges instead of on the
vertices, the model turns into the contact process with random edge
weights, which is first introduced by Chen and Yao in
\cite{Yao2012}, where a complete convergence theorem is shown.

\section{Main results}\label{section two}
In this section we give our main results. First we introduce some
notations and definitions. We assume that $\{\rho(x)\}_{x\in
\mathbb{Z}_+^d}$ are defined under the probability space
$(\Omega_d,\mathcal{F}_d,\mu_d)$. The expectation with respect to
$\mu_d$ is denoted by $E_{\mu_d}$. For $\omega\in \Omega_d$, we
denote by $P_{\lambda,\omega}$ the probability measure of our model
with vertex weights $\{\rho(x,\omega)\}_{x\in \mathbb{Z}_+^d}$.
$P_{\lambda,\omega}$ is called the quenched measure. The expectation
with respect to $P_{\lambda,\omega}$ is denoted by
$E_{\lambda,\omega}$. We define
\[
P_{\lambda,d}(\cdot)=E_{\mu_d}\big[P_{\lambda,\omega}(\cdot)\big]=\int
P_{\lambda,\omega}(\cdot)~\mu_d(d\omega),
\]
which is called the annealed measure. The expectation with respect
to $P_{\lambda,d}$ is denoted by $E_{\lambda,d}$.

For any $A\subseteq \mathbb{Z}_+^d$, we write $C_t$ as $C_t^A$ when
$C_0=A$. If $A=\{x\}$ for some $x\in \mathbb{Z}_+^d$, we write
$C_t^A$ as $C_t^x$ instead of $C_t^{\{x\}}$.

For $\lambda>\frac{1}{E(\rho^2)}$, where $E$ is the expectation with
respect to $\rho$, there is a unique solution $\theta>0$ to the
equation
\begin{equation}\label{equ 2.1}
E\big(\frac{\lambda\rho^2}{1+\lambda\rho\theta}\big)=1.
\end{equation}
Now we give the main result of this paper.
\begin{theorem}\label{theorem 2.1}
For any $\lambda>\frac{1}{E(\rho^2)}$ and $\theta$ defined as in
Equation \eqref{equ 2.1},
\[
\lim_{d\rightarrow+\infty}P_{\lambda,d}\big(C_t^O\neq
\emptyset,\forall~t\geq
0\big)=E\big(\frac{\lambda\rho\theta}{1+\lambda\rho\theta}\big).
\]
\end{theorem}

Theorem \ref{theorem 2.1} gives the asymptotic behavior of the
survival probability of the process conditioned on $O$ is the unique
initially infected vertex as the dimension $d$ grows to infinity.
The theorem only deals with the case where
$\lambda>\frac{1}{E(\rho^2)}$ because
\[
\lim_{d\rightarrow+\infty}P_{\lambda,d}\big(C_t^O\neq
\emptyset,\forall~t\geq 0\big)=0
\]
for any $\lambda<\frac{1}{E(\rho^2)}$ according to the main theorem
given in \cite{Xue2015}, which shows that the critical value of the
infection rate of the model converges to $\frac{1}{E(\rho^2)}$ as
$d\rightarrow+\infty$.

When $\rho\equiv1$, we have the following direct corollary.
\begin{corollary}\label{corollary 2.2}
If $\rho\equiv1$ and $\lambda>1$, then
\[
\lim_{d\rightarrow+\infty}P_{\lambda,d}\big(C_t^O\neq
\emptyset,\forall~t\geq 0\big)=\frac{\lambda-1}{\lambda}.
\]
\end{corollary}

The counterpart of Corollary \ref{corollary 2.2} for the classic
contact process on the lattice is given in \cite{Xue2017}.

The counterpart of Theorem \ref{theorem 2.1} for the contact process
with random edges weights on the (un-oriented) lattice is given in
\cite{Xue2017b}. It is claimed in \cite{Xue2017b} that
\[
\lim_{d\rightarrow+\infty}P_{\lambda,d}\big(C_t^O\neq
\emptyset,\forall~t\geq 0\big)=\frac{\lambda E\rho-1}{\lambda E\rho}
\]
for the process with edge weights which are independent copies of
$\rho$ and infection rate $\lambda>\frac{1}{E\rho}$.

As a auxiliary tool for the proof of Theorem \ref{theorem 2.1}, we
introduce a SIR (susceptible-infected-recovered) model with random
vertex weights on $\mathbb{Z}_+^d$.

After the vertex weights $\{\rho(x)\}_{x\in \mathbb{Z}_+^d}$ are
given, the SIR model $\{(S_t,I_t)\}_{t\geq 0}$ is a continuous-time
Markov process with state spapce
\[
X_2=\big\{(S,I):~S,I\subseteq \mathbb{Z}_+^d,~S\bigcap
I=\emptyset\big\}
\]
and transition rates function given by
\begin{equation}\label{equ 2.2 transition rate for SIR}
(S_t,I_t)\rightarrow
\begin{cases}
(S_t,I_t\setminus\{x\}) & \text{~at rate~} 1 \text{~if~} x\in I_t,\\
(S_t\setminus\{x\},I_t\bigcup\{x\}) & \text{~at rate~}
\frac{\lambda}{d}\sum\limits_{y:y\rightarrow
x}\rho(x)\rho(y)1_{\{y\in I_t\}} \text{~if~}x\in S_t.
\end{cases}
\end{equation}

For the SIR model, an infected vertex waits for an exponential time
with rate one to become recovered while a recovered vertex can never
be infected again.

We write $(S_t,I_t)$ as $(S_t^A,I_t^A)$ when
$(S_0,I_0)=(\mathbb{Z}_+^d\setminus A,A)$, then it is easy to check
that
\[
P_{\lambda,d}\big(I_t^O\neq \emptyset,\forall~t\geq 0\big)\leq
P_{\lambda,d}\big(C_t^O\neq \emptyset,\forall~t\geq 0\big).
\]
One way to check this inequality is to utilize the basic coupling of
Markov processes (see Section 3.1 of \cite{Lig1985}), we omit the
details. As a result, to prove Theorem \ref{theorem 2.1}, we only
need to show that
\begin{equation}\label{equ 2.3}
\liminf_{d\rightarrow+\infty}P_{\lambda,d}\big(I_t^O\neq
\emptyset,\forall~t\geq 0\big) \geq
E\big(\frac{\lambda\rho\theta}{1+\lambda\rho\theta}\big)
\end{equation}
and
\begin{equation}\label{equ 2.4}
\limsup_{d\rightarrow+\infty}P_{\lambda,d}\big(C_t^O\neq
\emptyset,\forall~t\geq 0\big) \leq
E\big(\frac{\lambda\rho\theta}{1+\lambda\rho\theta}\big).
\end{equation}

The proof of Theorem \ref{theorem 2.1} is divided into three
sections. In Section \ref{section three}, we introduce a branching
process $\{W_n\}_{n\geq 0}$ with random vertex weights on the
oriented rooted tree $\mathbb{T}^d$. We will show that the
probability that the process survives converges to
$E\big(\frac{\lambda\rho\theta}{1+\lambda\rho\theta}\big)$ as
$d\rightarrow+\infty$.

In Section \ref{section four}, we give the proof of Equation
\eqref{equ 2.3}. The proof relies on a coupling relationship between
the branching process and the SIR model. A technique introduced in
\cite{Xue2015} is utilized.

In Section \ref{section five}, we give the proof of Equation
\eqref{equ 2.4}. The proof relies on a coupling relationship of the
three aforesaid processes.

\section{A branching process with vertex weights}\label{section three}
In this section we introduce a branching process with random vertex
weights on the oriented rooted tree. We denote by $\mathbb{T}^d$ the
rooted tree that the root has $d$ neighbors while any other vertex
on the tree has $d+1$ neighbors. We denote by $\Upsilon$ the root of
the tree. There is a function $f:\mathbb{T}^d\rightarrow
\{0,1,2,\ldots\}$ satisfies the following conditions.

(1) $f(\Upsilon)=0$.

(2) $f(x)=1$ for each neighbor $x$ of $\Upsilon$.

(3) For any $y\neq \Upsilon$, there is one neighbor $u$ of $y$ that
$f(u)=f(y)-1$ while there are $d$ neighbors $v$ of $y$ that
$f(v)=f(y)+1$.

For $x,y\in \mathbb{T}^d$, we write $x\Rightarrow y$ when and only
when $x$ and $y$ are neighbors and $f(y)=f(x)+1$.

Intuitively, $\Upsilon$ is the ancestor of a family and has $d$
sons. Each other individual in this family has one father and $d$
sons. $x\Rightarrow y$ when and only when $y$ is a son of $x$.

We assume that $\{\rho(x)\}_{x\in \mathbb{T}^d}$ are i.i.d. copies
of the random variable $\rho$, which is defined as in Section
\ref{section one}. After the vertex weights are given, we assume
that $Y(x)$ is an exponential time with rate one for each $x\in
\mathbb{T}^d$ while $U(x,y)$ is an exponential time with rate
$\frac{\lambda}{d}\rho(x)\rho(y)$ for any $x,y\in \mathbb{T}^d$ that
$x\Rightarrow y$. We assume that all these exponential times are
independent under the given vertex weights. Then, the branching
process $\{W_n\}_{n\geq 0}$ is defined as follows.

(1) $W_0=\Upsilon$.

(2) For $n\geq 0$, $W_{n+1}=\big\{y:x\Rightarrow y
\text{~and~}U(x,y)<Y(x)\text{~for some~}x\in W_n\big\}$.

$\{W_n\}_{n\geq 0}$ describes the spread of a SIR epidemic on
$\mathbb{T}^d$. Initially, $\Upsilon$ is infected. A healthy vertex
may only be infected by its father. If $x$ is infected, then $x$
waits for an exponential time with rate one to become recovered
while waits for an exponential time with rate
$\frac{\lambda}{d}\rho(x)\rho(y)$ to infect the son $y$. The
infection really occurs when and only when $y$ is infected before
the moment $x$ is recovered, i.e., $U(x,y)<Y(x)$.

Similar with what we have done in Section \ref{section two}, we
denote by $\widehat{P}_{\lambda,\omega}$ the quenched measure of the
branching process with respect to the random environment $\omega$ in
the space where $\{\rho(x)\}_{x\in \mathbb{T}^d}$ are defined. We
denote by $\widehat{P}_{\lambda,d}$ the annealed measure. Note that
according to our definition, for $x\Rightarrow y\Rightarrow z$,
$U(x,y)$ and $U(y,z)$ are independent under
$\widehat{P}_{\lambda,\omega}$ while positively correlated under
$\widehat{P}_{\lambda,d}$.

The branching process $\{W_n\}_{n\geq 0}$ with random vertex weights
on the oriented tree $\mathbb{T}^d$ is first introduced in
\cite{Pan2017}. Some results obtained in \cite{Pan2017} will be
directly utilized in this section.

The following lemma is crucial for us to prove Theorem \ref{theorem
2.1}.
\begin{lemma}\label{lemma 3.1}
For any $\lambda>\frac{1}{E(\rho^2)}$ and $\theta$ defined as in
Equation \eqref{equ 2.1},
\[
\lim_{d\rightarrow+\infty}\widehat{P}_{\lambda,d}\big(W_n\neq
\emptyset,\forall ~n\geq
0\big)=E\big(\frac{\lambda\rho\theta}{1+\lambda\rho\theta}\big).
\]
\end{lemma}
The remainder of this section is devoted to the proof of Lemma
\ref{lemma 3.1}. From now on we assume that
$\lambda>\frac{1}{E(\rho^2)}$. For any $s\in [0.M]$, we define
\[
F_d(s)=\widehat{P}_{\lambda,d}\Big(W_n=\emptyset\text{~for
some~}n\geq 0\Big|\rho(\Upsilon)=s\Big),
\]
then the following two lemmas are crucial for us to prove Lemma
\ref{lemma 3.1}.
\begin{lemma}\label{lemma 3.2}
If $\{d_l\}_{l\geq 1}$ is a subsequence of $1,2,3,\ldots$ such that
\[
\lim_{l\rightarrow+\infty}F_{d_l}(s)\exists:=F(s)
\]
for any $s\in [0,M]$, then
\[
F(s)=\frac{1}{1+\lambda s\theta}.
\]
\end{lemma}

\begin{lemma}\label{lemma 3.3}
For $d\geq 1$ and $0\leq s<t\leq M$,
\[
|F_d(s)-F_d(t)|\leq \lambda (t-s)M.
\]
\end{lemma}

We first show how to utilize Lemmas \ref{lemma 3.2} and \ref{lemma
3.3} to prove Lemma \ref{lemma 3.1}. The proofs of Lemmas \ref{lemma
3.2} and \ref{lemma 3.3} are given at the end of this section.

\proof[Proof of Lemma \ref{lemma 3.1}]

If Lemma \ref{lemma 3.1} does not hold, then there are a constant
$\epsilon_0>0$ and a subsequence $\{a_l\}_{l\geq 1}$ of
$1,2,3,\ldots$ that
\begin{equation}\label{equ 3.4}
|\widehat{E}_{\lambda,a_l}\big(F_{a_l}(\rho)\big)-E\big(\frac{1}{1+\lambda\rho\theta}\big)|>\epsilon,
\end{equation}
since
\[
\widehat{P}_{\lambda,d}\big(W_n\neq \emptyset,\forall ~n\geq
0\big)=1-E_{\lambda,d}\big(F_d(\rho)\big).
\]
Since $0\leq F_d(\cdot)\leq 1$, according to a classic procedure of
picking subsequences, there is a subsequence $\{d_j\}_{j\geq 1}$ of
$\{a_l\}_{l\geq 1}$ such that
\[
\lim_{j\rightarrow+\infty}F_{d_j}(r)\exists:=F_\Delta(r)
\]
for any $r\in \mathbb{Q}$. It is obviously that $F_d(s)$ is
decreasing with $s$ for each $d\geq 1$, then
\[
F_d(r_1)\geq F_d(s)\geq F_d(r_2)\text{~and~}F_\Delta(r_1)\geq
F_\Delta(r_2)
\]
for any $r_1<s<r_2,r_1,r_2\in \mathbb{Q}$. As a result, it is
reasonable to define
\[
F_\Delta^-(s)=\lim_{r\uparrow s,r\in
\mathbb{Q}}F_\Delta(r)\text{~and~} F_\Delta^+(s)=\lim_{r\downarrow
s,r\in \mathbb{Q}}F_\Delta(r)
\]
for any $s\not \in \mathbb{Q}$ and hence
\[
\limsup_{j\rightarrow+\infty}F_{d_j}(s)\leq
F_\Delta^-(s)\text{~while~}\liminf_{j\rightarrow+\infty}F_{d_j}(s)\geq
F_\Delta^-(s).
\]
By Lemma \ref{lemma 3.3},
\[
|F_\Delta(r_1)-F_\Delta(r_2)|\leq \lambda M(r_2-r_1)
\]
for $r_1<s<r_2,r_1,r_2\in \mathbb{Q}$. Therefore, let $r_1\uparrow
s$ and $r_2\downarrow s$,
\[
F_\Delta^-(s)=F_\Delta^+(s):=F_\Delta(s)
\]
and
\[
\lim_{j\rightarrow+\infty}F_{d_j}(s)=F_\Delta(s)
\]
for any $s\not \in \mathbb{Q}$. As a result,
\[
\lim_{j\rightarrow+\infty}F_{d_j}(s)\exists:=F_\Delta(s)
\]
for any $s\in [0,M]$. Then, by Lemma \ref{lemma 3.2},
\[
F_\Delta(s)=\frac{1}{1+\lambda s\theta}
\]
for any $s\in [0,M]$ and hence
\[
\lim_{j\rightarrow+\infty}\widetilde{E}_{\lambda,d_j}\big(F_{d_j}(\rho)\big)=E\big(\frac{1}{1+\lambda\rho\theta}\big).
\]
However, this is contradictory with Equation \eqref{equ 3.4} since
$\{d_j\}_{j\geq 1}$ is a subsequence of $\{a_l\}_{l\geq 1}$. As a
result, Lemma \ref{lemma 3.1} holds and the proof is complete.

\qed

At last we give the proof of Lemmas \ref{lemma 3.2} and \ref{lemma
3.3}.

\proof[Proof of Lemma \ref{lemma 3.2}]

For $\Upsilon\Rightarrow y$, conditioned on $Y(\Upsilon),
\rho(\Upsilon),\rho(y)$, the probability that $\Upsilon$ infects $y$
is
\begin{align*}
&P_{\lambda,\omega}\big(U(\Upsilon,y)<Y(\Upsilon)\big|\rho(\Upsilon),Y(\Upsilon),
\rho(y)\big) \\
&=1-e^{-\frac{\lambda}{d}\rho(\Upsilon)\rho(y)Y(\Upsilon)}.
\end{align*}
If $\{W_n\}_{n\geq 0}$ dies out, then for any $y$ such that
$\Upsilon$ infects $y$, the epidemic on the subtree consisted of $y$
and its descendant must die out, the probability of which is
$F_d(\rho(y))$. As a result,
\begin{align*}
&\widehat{E}_{\lambda,d}\Big[W_n=\emptyset\text{~for some~}n\geq 0\Big|\rho(\Upsilon),Y(\Upsilon),\{\rho(y):\Upsilon\Rightarrow y\}\Big]\\
&=\prod_{y:\Upsilon\Rightarrow
y}\Big(F_d(\rho(y))(1-e^{-\frac{\lambda}{d}\rho(\Upsilon)\rho(y)Y(\Upsilon)})+e^{-\frac{\lambda}{d}\rho(\Upsilon)\rho(y)Y(\Upsilon)}\Big)
\end{align*}
and hence
\[
F_d(s)=\widehat{E}_{\lambda,d}\Bigg[\prod_{y:\Upsilon\Rightarrow
y}\Big(F_d(\rho(y))(1-e^{-\frac{\lambda}{d}s\rho(y)Y(\Upsilon)})+e^{-\frac{\lambda}{d}s\rho(y)Y(\Upsilon)}\Big)\Bigg].
\]
Since $\{\rho(y):\Upsilon\Rightarrow y\}$ are independent,
\begin{equation}\label{equ 3.1}
F_d(s)=\widehat{E}_{\lambda,d}\Bigg[\Big(H_d(Y(\Upsilon))\Big)^d\Bigg]=E\Bigg[\Big(H_d(Y_0)\Big)^d\Bigg],
\end{equation}
where
\[
H_d(t)=E\Big(F_d(\rho)(1-e^{-\frac{\lambda}{d}st\rho})+e^{-\frac{\lambda}{d}st\rho}\Big)
\]
for any $t\geq 0$ and $Y_0$ is an exponential time with rate one
defined under some space we do not care. According to our
assumption, it is easy to check that
\begin{equation}\label{equ 3.2}
\lim_{l\rightarrow+\infty}d_l(H_{d_l}(t)-1)=-\lambda
stE\Big(\rho(1-F(\rho))\Big).
\end{equation}
According to the theory of calculus, if $a_d\rightarrow 0,
c_d\rightarrow +\infty$ and $a_dc_d\rightarrow c$, then
$(1+a_d)^{c_d}\rightarrow e^c$. Therefore, by Equations \eqref{equ
3.1} and \eqref{equ 3.2},
\begin{equation}\label{equ 3.3}
\lim_{l\rightarrow+\infty}F_{d_l}(s)=Ee^{-\lambda
sY_0\widetilde{\theta}}=\frac{1}{1+\lambda s\widetilde{\theta}},
\end{equation}
where $\widetilde{\theta}=E\big(\rho(1-F(\rho))\big)$. As a result,
$F(s)=\frac{1}{1+\lambda s \widetilde{\theta}}$ for any $s$ and we
only need to show that $\widetilde{\theta}=\theta$. According to the
definition of $\widetilde{\theta}$,
\[
\widetilde{\theta}=E\Big(\rho(1-F(\rho))\Big)=E\Big(\rho(1-\frac{1}{1+\lambda
\rho\widetilde{\theta}})\Big)=E\big(\frac{\lambda\rho^2\widetilde{\theta}}{1+\lambda\rho
\widetilde{\theta}}\big).
\]
Therefore, to prove $\theta=\widetilde{\theta}$ we only need to show
that $\widetilde{\theta}\neq 0$. This fact follows directly from the
conclusion that
\[
\limsup_{d\rightarrow+\infty}E\big(F_d(\rho)\big)<1
\]
when $\lambda>\frac{1}{E(\rho^2)}$, which is proved in
\cite{Pan2017}.

\qed

\proof[Proof of Lemma \ref{lemma 3.3}]

We denote by $\{W_n^s\}_{n\geq 0}$ the branching process conditioned
on $\rho(\Upsilon)=s$ and denote by $\{W_n^t\}_{n\geq 0}$ the
branching process conditioned on $\rho(\Upsilon)=t$. We couple these
two branching processes in a same probability space as follows. For
any $x\in \mathbb{T}^d$, we assume that these two processes utilize
the same exponential time $Y(x)$ with rate one. For any $x\neq
\Upsilon$ and $x\Rightarrow z$, we assume that these two processes
utilize the same exponential time $U(x,z)$ with rate
$\frac{\lambda}{d}\rho(x)\rho(z)$. For each $y$ that
$\Upsilon\Rightarrow y$, we assume that $\{W_n^s\}_{s\geq 0}$
utilizes an exponential time $U_s(\Upsilon,y)$ with rate
$\frac{\lambda}{d}s\rho(y)$ while $\{W_n^t\}_{t\geq 0}$ utilizes an
exponential time
\[
U_t(\Upsilon,y)=\inf\big\{U_s(\Upsilon,y),U_{t-s}(\Upsilon,y)\big\},
\]
where $U_{t-s}(\Upsilon,y)$ is an exponential time with rate
$\frac{\lambda}{d}(t-s)\rho(y)$ and is independent of
$U_s(\Upsilon,y), Y(\Upsilon)$ under the quenched measure.
Therefore, $U_t(\Upsilon,y)$ is an exponential time with rate
$\frac{\lambda}{d}t\rho(y)$. According to the coupling of
$\{W_n^s\}_{n\geq 0}$ and $\{W_n^t\}_{n\geq 0}$,
\begin{align*}
|F_d(t)-F_d(s)|&=P\big(\{W_n^t\}_{n\geq 0}\text{~survives while~}\{W_n^s\}_{n\geq 0}\text{~dies out}\big)\\
&\leq P(W_1^s\neq W_1^t) \\
&=P\big(U_{t-s}(\Upsilon,y)<Y(\Upsilon)<U_s(\Upsilon,y)\text{~for some~}y\big)\\
&\leq \sum_{y:\Upsilon\Rightarrow
y}P\big(U_{t-s}(\Upsilon,y)<Y(\Upsilon)<U_s(\Upsilon,y)\big)\\
&=\sum_{y:\Upsilon\Rightarrow y}\widehat{E}_{\lambda,d}\Big[e^{-\frac{\lambda}{d}s\rho(y)Y(\Upsilon)}-e^{-\frac{\lambda}{d}t\rho(y)Y(\Upsilon)}\Big]\\
&=\sum_{y:\Upsilon\Rightarrow
y}\widehat{E}_{\lambda,d}\Big[\frac{1}{1+\frac{\lambda}{d}s\rho(y)}-\frac{1}{1+\frac{\lambda}{d}t\rho(y)}\Big]\\
&=dE\Big[\frac{\frac{\lambda}{d}(t-s)\rho}{(1+\frac{\lambda
s\rho}{d})(1+\frac{\lambda t\rho}{d})}\Big]\leq
d\frac{\lambda}{d}(t-s)M=\lambda(t-s)M
\end{align*}
and the proof is complete.

\qed

\section{Proof of Equation \eqref{equ 2.3}}\label{section four}
In this section we give the proof of Equation \eqref{equ 2.3}. For
later use, we assume that there exists $\epsilon>0$ that
\begin{equation}\label{equ 4.1}
P\big(\rho=0\text{~or~}\rho\in [\epsilon,M]\big)=1,
\end{equation}
where $M$ is defined as in Section \ref{section one}.

This assumption is without loss of generality according to the
following analysis. For $\rho$ not satisfying \eqref{equ 4.1}, we
let $\rho_m=\rho1_{\{\rho\geq 1/m\}}$, then $\rho\geq \rho_m$ and
$\rho_m\rightarrow \rho$ as $m\rightarrow+\infty$. It is obviously
that $\rho_m$ satisfies \eqref{equ 4.1} while the process with
weights respect to $\rho$ has larger probability to survive than
that with weights respect to $\rho_m$. As a result, if Equation
\eqref{equ 2.3} holds under assumption \eqref{equ 4.1}, then
\[
\liminf_{d\rightarrow+\infty} P_{\lambda,d,\rho}\big(I_t^O\neq
\emptyset,\forall~t\geq 0\big)\geq \liminf_{d\rightarrow+\infty}
P_{\lambda,d,\rho_m}\big(I_t^O\neq \emptyset,\forall~t\geq
0\big)\geq
E\big(\frac{\lambda\rho_m\theta_m}{1+\lambda\rho_m\theta_m}\big)
\]
for any sufficiently large $m$, where $P_{\lambda,d,\rho}$ is the
annealed measure of the process with vertex weights which are i.i.d
copies of $\rho$ while $\theta_m$ satisfies
\[
E\big(\frac{\lambda\rho_m^2}{1+\lambda\rho_m\theta_m}\big)=1
\]
and it is easy to check that
$\lim_{m\rightarrow+\infty}\theta_m=\theta$. Let $m\rightarrow
+\infty$, then Equation \eqref{equ 2.3} holds for general $\rho$.

First we give a sketch of the proof, which is inspired by the
approach introduced in \cite{Xue2017b}. We divide $\mathbb{Z}_+^d$
into two parts $\Gamma_1$ and $\Gamma_2$ such that $\Gamma_1\bigcap
\Gamma_2=\emptyset$. The first step is to show that with probability
at least
$E\big(\frac{\lambda\rho\theta}{1+\lambda\rho\theta}\big)+o(1)$
there exists a path starting at $O$ on $\Gamma_1$ with length
$O(\log d)$ that all the vertices on this path have ever been
infected. The second step is to show that conditioned on the
existence of the aforesaid path, the vertices on this path infects
$K(d)$ vertices on $\Gamma_2$ through edges connecting $\Gamma_1$
and $\Gamma_2$ with high probability, where $K(d)\rightarrow+\infty$
as $d\rightarrow+\infty$. The third step is to show that conditioned
on $K(d)$ vertices are initially infected on $\Gamma_2$, the SIR
model on $\Gamma_2$ survives with high probability. To prove the
first step, we construct a couple between the SIR on $Z_+^d$ and the
branching process introduced in Section \ref{section three}.

To give our proof, we introduce some definitions and notations. For
sufficiently large $d$, we define $N(d)=\log(\log d)$ and
\begin{align*}
&\Gamma_1=\Big\{x=(x_1,x_2,\ldots,x_d)\in \mathbb{Z}_+^d:\sum_{i=d-\lfloor\frac{d}{N(d)}\rfloor+1}^dx_i=0\Big\},\\
&\Gamma_2=\Big\{x\in
\mathbb{Z}_+^d:\sum_{i=d-\lfloor\frac{d}{N(d)}\rfloor+1}^dx_i>0\Big\},\\
&\Gamma_3=\Big\{x\in
\mathbb{Z}_+^d:\sum_{i=d-\lfloor\frac{d}{N(d)}\rfloor+1}^dx_i=1\Big\}\subseteq
\Gamma_2.
\end{align*}
For $n\geq 0$, we define
\[
V_n=\Big\{x\in \mathbb{Z}_+^d:\|x\|=n\text{~and~}x\in
I_t^O\text{~for some~}t\geq 0\Big\}
\]
as the set of vertices which have ever been infected with $l_1$ norm
$n$. Since in the SIR model, infection can not occur repeatedly
between neighbors, $\{V_n\}_{n\geq 0}$ can be defined equivalently
as the following way. For each $x\in \mathbb{Z}^d_+$, let
$\widetilde{Y}(x)$ be an exponential time with rate one. For any
$x,y$ that $x\rightarrow y$, let $\widetilde{U}(x,y)$ be an
exponential time with rate $\frac{\lambda}{d}\rho(x)\rho(y)$. We
assume that all these exponential times are independent under the
quenched measure with respect to the given edge weights, then

(1) $V_0=\{O\}$.

(2) For each $n\geq 0$,
\[
V_{n+1}=\big\{y\in \mathbb{Z}_+^d:x\rightarrow y
\text{~and~}\widetilde{U}(x,y)<\widetilde{Y}(x)\text{~for some~}x\in
V_n\big\}.
\]

The intuitive explanation of the above definition is similar with
that of the branching process introduced in Section \ref{section
three}. $\widetilde{Y}(x)$ is time $x$ waits for to become recovered
after $x$ is infected while $\widetilde{U}(x,y)$ is the time $x$
waits for to infect $y$.

For later use, we define $\{\widehat{V}_n\}_{n\geq 0}$ as the
vertices which have ever been infected for the SIR model confined on
$\Gamma_1$. In details,

(1) $\widehat{V}_0=\{O\}$.

(2) For each $n\geq 0$,
\[
\widehat{V}_{n+1}=\big\{y\in \Gamma_1:x\rightarrow y
\text{~and~}\widetilde{U}(x,y)<\widetilde{Y}(x)\text{~for some~}x\in
\widehat{V}_n\big\}.
\]

Let $\sigma>0$ be an arbitrary constant that $\sigma\log(\lambda
M^2)<\frac{1}{10}$, then we have the following lemma.

\begin{lemma}\label{lemma 4.1}
For $\lambda>\frac{1}{E(\rho^2)}$,
\[
\liminf_{d\rightarrow+\infty}P_{\lambda,d}\Big(\widehat{V}_{\lfloor\sigma\log
d\rfloor}\neq \emptyset\Big)\geq
E\big(\frac{\lambda\rho\theta}{1+\lambda\rho\theta}\big).
\]
\end{lemma}

The proof of Lemma \ref{lemma 4.1} is given in Subsection
\ref{subsection 4.2}. As a preparation of this proof, we give a
coupling of $\{W_n\}_{n\geq 0}$ and $\{V_n\}_{n\geq 0}$ in
Subsection \ref{subsection 4.1}.

To execute the second step as we have introduced, we define
\[
D=\big\{y\in \Gamma_3: x\rightarrow
y\text{~and~}\widetilde{U}(x,y)<\widetilde Y(x)\text{~for some~}x\in
\bigcup_{n=0}^{\lfloor\sigma\log d\rfloor}\widehat{V}_n\big\}
\]
as the set of vertices in $\Gamma_3$ which have been infected by
vertices in $ \bigcup_{n=0}^{\lfloor\sigma\log
d\rfloor}\widehat{V}_n$ through edges connecting $\Gamma_1$ and
$\Gamma_2$. Let $K(d)=\lfloor\frac{\sqrt{\log d}}{N(d)}\rfloor$,
then we have the following lemma.
\begin{lemma}\label{lemma 4.2}
For $\lambda>\frac{1}{E(\rho^2)}$,
\[
\lim_{d\rightarrow+\infty}P_{\lambda,d}\Big(|D|\geq
K(d)\Big|\widehat{V}_{\lfloor\sigma\log d\rfloor}\neq
\emptyset\Big)=1,
\]
where $|D|$ is the cardinality of $D$.
\end{lemma}
The proof of Lemma \ref{lemma 4.2} is given in Subsection
\ref{subsection 4.3}.

As we have introduced, the third step of the proof of Equation \eqref{equ 2.3} is to show that the SIR model survives with high probability conditioned on $K(d)$ vertices are initially infected. In detail, we have the following
lemma.
\begin{lemma}\label{lemma 4.3}
If $\lambda>\frac{1}{E(\rho^2)}$, then there exists $m(d)$ for each
$d\geq 1$ such that
\[
\lim_{d\rightarrow+\infty}m(d)=1
\]
and
\[
P_{\lambda,d}\Big(I_t^A\neq \emptyset, \forall~t\geq
0\Big)\geq m(d)
\]
for any $A\subseteq \Gamma_3$ with $|A|=K(d)$.
\end{lemma}
The proof of Lemma \ref{lemma 4.3} is given in Subsection
\ref{subsection 4.4}.

Now we show how to utilize Lemmas \ref{lemma 4.1}, \ref{lemma 4.2}
and \ref{lemma 4.3} to prove Equation \eqref{equ 2.3}.

\proof[Proof of Equation \eqref{equ 2.3}]

For $x,y\in \mathbb{Z}_+^d$, we write $x\rightrightarrows y$ when
there exists $x_1,x_2,\ldots,x_m$ for some integer $m\geq 1$ such
that $x=x_0\rightarrow x_1\rightarrow
x_2\rightarrow\ldots\rightarrow x_m\rightarrow x_{m+1}=y$ and
$\widetilde{U}(x_j,x_{j+1})<\widetilde{Y}(x_j)$ for all $0\leq j\leq
m$. Then, according to the meanings of the exponential times
$\widetilde{U}(\cdot,\cdot)$ and $\widetilde{Y}(\cdot)$,
\[
\bigcup_{t\geq 0}I_t^A=A\bigcup\big\{y:x\rightrightarrows
y\text{~for some~}x\in A\big\}.
\]
Since each infected vertex becomes recovered in an exponential time
with rate one, the infected vertices never die out when and only
when there are infinite many vertices have ever been infected.
Therefore,
\begin{equation}\label{equ 4.01}
\Big\{I_t^A\neq \emptyset,\forall~t\geq
0\Big\}=\Big\{|\{y:x\rightrightarrows y\text{~for some~}x\in
A\}|=+\infty\Big\}
\end{equation}
for any finite $A$. According to the definition of $D$,
$O\rightrightarrows x$ for any $x\in D$. As a result,
\begin{equation}\label{equ 4.02}
\Big\{y:x\rightrightarrows y \text{~for some~}x\in D\Big\}\subseteq
\Big\{y:O\rightrightarrows y\Big\}\subseteq \bigcup_{t\geq 0}I_t^O,
\end{equation}
since $O\rightrightarrows y$ when $O\rightrightarrows x$ and
$x\rightrightarrows y$. By Equations \eqref{equ 4.01} and \eqref{equ
4.02},
\begin{equation}\label{equ 4.03}
P_{\lambda,d}\Big(|\{y:x\rightrightarrows y\text{~for some~}x\in
D\}|=+\infty\Big)\leq P_{\lambda,d}\big(I_t^O\neq
\emptyset,\forall~t\geq 0\big).
\end{equation}
According to the conditional probability formula,
\begin{align}\label{equ 4.04}
&P_{\lambda,d}\Big(|\{y:x\rightrightarrows y\text{~for some~}x\in
D\}|=+\infty\Big)\\
&\geq  P_{\lambda,d}\Big(|\{y:x\rightrightarrows y\text{~for
some~}x\in D\}|=+\infty\Big||D|\geq
K(d),\widehat{V}_{\lfloor\sigma\log d\rfloor}\neq \emptyset\Big)\notag\\
&\times P_{\lambda,d}\Big(|D|\geq
K(d)\Big|\widehat{V}_{\lfloor\sigma\log d\rfloor}\neq
\emptyset\Big)P_{\lambda,d}\Big(\widehat{V}_{\lfloor\sigma\log
d\rfloor}\neq \emptyset\Big).\notag
\end{align}
Note that $D$ is a subset of $\Gamma_3$ with cardinality at least $K(d)$ while the event $\big\{|\{y:x\rightrightarrows y\text{~for
some~}x\in D\}|=+\infty\big\}$ is positively correlated with the event $\big\{|D|\geq
K(d),\widehat{V}_{\lfloor\sigma\log d\rfloor}\neq \emptyset\big\}$. As a result, by
Lemma \ref{lemma 4.3} and Equation \eqref{equ 4.04},
\begin{align}\label{equ 4.05}
&P_{\lambda,d}\Big(|\{y:x\rightrightarrows y\text{~for some~}x\in
D\}|=+\infty\Big||D|\geq K(d),\widehat{V}_{\lfloor\sigma\log
d\rfloor}\neq \emptyset\Big)\notag\\
&\geq \inf\Bigg\{P_{\lambda,d}\Big(|\{y:x\rightrightarrows
y\text{~for some~}x\in A\}|=+\infty\Big):\\
&|A|=K(d)\text{~and~}A\subseteq \Gamma_3\Bigg\}\notag\\
&=\inf\Bigg\{P_{\lambda,d}\Big(I_t^A\neq \emptyset,\forall~t\geq
0\Big):|A|=K(d)\text{~and~}A\subseteq
\Gamma_3\Bigg\}\geq m(d). \notag
\end{align}
By Equations \eqref{equ 4.04}, \eqref{equ 4.05}, Lemmas \ref{lemma
4.1}, \ref{lemma 4.2} and \ref{lemma 4.3},
\begin{align}\label{equ 4.06}
&\liminf_{d\rightarrow+\infty}P_{\lambda,d}\Big(|\{y:x\rightrightarrows
y\text{~for some~}x\in D\}|=+\infty\Big)\\
&\geq E\big(\frac{\lambda\rho\theta}{1+\lambda\rho\theta}\big)
\lim_{d\rightarrow+\infty}m(d)=E\big(\frac{\lambda\rho\theta}{1+\lambda\rho\theta}\big).\notag
\end{align}
Equation \eqref{equ 2.3} follows directly from Equations \eqref{equ
4.03} and \eqref{equ 4.06}.

\qed

\subsection{The coupling between $\{W_n\}_{n\geq 1}$ and $\{V_n\}_{n\geq 1}$}\label{subsection 4.1}
In this section, we give a coupling between the SIR model
$\{V_n\}_{n\geq 0}$ on $\mathbb{Z}^d_+$ and the branching process
$\{W_n\}_{n\geq 0}$ on $\mathbb{T}^d$.

We let $\{\rho(x)\}_{x\in \mathbb{Z}^d_+}$ be i.i.d copies of $\rho$
as defined in Section \ref{section one}. We let
$\{\widetilde{Y}(x)\}_{x\in \mathbb{Z}^d_+}$ and
$\{\widetilde{U}(x,y)\}_{x\rightarrow y}$ be exponential times with
respect to $\{\rho(x)\}_{x\in \mathbb{Z}^d_+}$ as defined at the
beginning of this section. We let $\{V_n\}_{n\geq 0}$ be the SIR
model with respect to $\widetilde{Y}(\cdot)$ and
$\widetilde{U}(\cdot,\cdot)$ as defined at the beginning of this
section. Now we give the evolution of $\{W_n\}_{n\geq 0}$ by
induction.

We let $W_0=\Upsilon$, $\rho(\Upsilon)=\rho(O)$ and
$Y(\Upsilon)=\widetilde{Y}(O)$, where $O$ is the origin of
$\mathbb{Z}_+^d$. For the $d$ sons denoted by $n_1,n_2,\ldots,n_d$
of $\Upsilon$, we let $\rho(n_i)=\rho(e_i),
Y(n_i)=\widetilde{Y}(e_i)$ and
$U(\Upsilon,n_i)=\widetilde{U}(O,e_i)$ for each $1\leq i\leq d$,
where $e_i$ is the elementary unit vector of $\mathbb{Z}_+^d$ as
defined in Section \ref{section one}.  Then $W_1$ is defined
according to the values of $\{U(\Upsilon,n_i)\}_{1\leq i\leq d}$ and
$Y(\Upsilon)$ as in Section \ref{section three}.

For $n\geq 1$, if $|V_n|=|W_n|$ and there is a bijection
$g_n:V_n\rightarrow W_n$ such that $\rho(g_n(x))=\rho(x)$ and
$Y(g_n(x))=\widetilde{Y}(x)$ for each $x\in V_n$, then we say that
our coupling is successful at step $n$. It is obviously that our
coupling is successful at step $n=1$ since $g_1$ can be defined as
$g_1(e_i)=n_i$ for any $e_i\in V_1$.

If $\{W_m\}_{m\leq n}$ is well defined and the coupling is
successful at step $m$ for all $1\leq m\leq n$, then $W_{n+1}$ is
defined as follows. For any $x\in V_n$, we define
\[
q(x)=\Big\{y:x\rightarrow y \text{~and~}z\rightarrow y\text{~for
some~}z\in V_n\setminus\{x\}\Big\},
\]
\[
\psi(x)=\big\{y:x\rightarrow y\big\}\setminus q(x)
\]
and $h(x)=d-|q(x)|$, then $|\psi(x)|=h(x)$. For each $x\in V_n$, we
arbitrarily choose $h(x)$ sons of $g_n(x)$, which are denoted by
$w_1, w_2,\ldots,w_{h(x)}$. Giving the $h(x)$ elements in $\psi(x)$
an arbitrary order $y_1,y_2,\ldots,y_{_{h(x)}}$, then we let
$\rho(w_i)=\rho(y_i)$, $Y(w_i)=\widetilde{Y}(y_i)$ and
$U(g_n(x),w_i)=\widetilde{U}(x,y_i)$ for each $1\leq i\leq h(x)$.
For any son $u$ of $g_n(x)$ which is not in
$\{w_1,w_2,\ldots,w_{h(x)}\}$, let $Y(u)$ be an exponential time
with rate one and $\rho(u)$ be an independent copy of $\rho$ that
$Y(u)$ and $\rho(u)$ are independent of the aforesaid exponential
times and vertex weights while let $U(g_n(x),u)$ be an exponential
time with rate $\frac{\lambda}{d}\rho(g_n(x))\rho(u)$. Then,
$W_{n+1}$ is defined according to the values of $\{Y(g_n(x))\}_{x\in
V_n}$ and $\{U(g_n(x),w)\}_{x\in V_n,g_n(x)\Rightarrow w}$ as in
Section \ref{section three}.

If $n\geq 2$ is the first step that the coupling is not successful,
then we let $\{W_m\}_{m\geq n+1}$ evolves independently of
$\{V_m\}_{m\geq n+1}$.

From now on we assume that $\{W_n\}_{n\geq 0}$ and $\{V_n\}_{n\geq
0}$ are defined under the same probability space. The annealed
measure is still denoted by $P_{\lambda,d}$.

The remainder of this subsection is devoted to the proof of the
following lemma.
\begin{lemma}\label{lemma 4.1.1}
For $\lambda>\frac{1}{E(\rho^2)}$ and $\sigma\in
(0,\frac{1}{10\log(\lambda M^2)})$,
\[
\liminf_{d\rightarrow+\infty}P_{\lambda,d}\big(V_{\lfloor\sigma\log
d\rfloor}\neq \emptyset\big)\geq
E\big(\frac{\lambda\rho\theta}{1+\lambda\rho\theta}\big).
\]
\end{lemma}

The proof of Lemma \ref{lemma 4.1.1} relies heavily on the following
lemma.

\begin{lemma}\label{lemma 4.1.2}
For given $\sigma\in (0,\frac{1}{10\log(\lambda M^2)})$, we denote
by $B(d)$ the event that the coupling of $\{V_n\}_{n\geq 0}$ on
$\mathbb{Z}_+^d$ and $\{W_n\}_{n\geq 0}$ on $\mathbb{T}^d$ is
successful at step $m$ for all
\[
1\leq m \leq \lfloor\sigma\log d\rfloor,
\]
then
\[
\lim_{d\rightarrow+\infty}P_{\lambda,d}\big(B(d)\big)=1.
\]
\end{lemma}

The proof of Lemma \ref{lemma 4.1.2} is given at the end of this
subsection. Now we show how to utilize Lemma \ref{lemma 4.1.2} to
prove Lemma \ref{lemma 4.1.1}.

\proof[Proof of Lemma \ref{lemma 4.1.1}]

On the event $B(d)$, $|W_m|=|V_m|$ for all $1\leq m\leq
\lfloor\sigma\log d\rfloor$. As a result,
\[
P_{\lambda,d}\big(V_{\lfloor\sigma\log d\rfloor}\neq \emptyset\big)
\geq P_{\lambda,d}\big(W_{\lfloor\sigma\log d\rfloor}\neq
\emptyset\big)-P_{\lambda,d}\big(B(d)^c\big),
\]
where $B(d)^c$ is the complementary set of $B(d)$. Therefore, by
Lemmas \ref{lemma 3.1} and \ref{lemma 4.1.2},
\begin{align*}
&\liminf_{d\rightarrow+\infty}
P_{\lambda,d}\big(V_{\lfloor\sigma\log d\rfloor}\neq
\emptyset\big)\\
&\geq \liminf_{d\rightarrow+\infty}
P_{\lambda,d}\big(W_{\lfloor\sigma\log d\rfloor}\neq
\emptyset\big)-\lim_{d\rightarrow+\infty}P_{\lambda,d}\big(B(d)^c\big)\\
&=\liminf_{d\rightarrow+\infty}
P_{\lambda,d}\big(W_{\lfloor\sigma\log d\rfloor}\neq
\emptyset\big)\\
&\geq \lim_{d\rightarrow+\infty}P_{\lambda,d}\big(W_n\neq
\emptyset,\forall~n\geq
0\big)=E\big(\frac{\lambda\rho\theta}{1+\lambda\rho\theta}\big)
\end{align*}
and the proof is complete.

\qed

At the end of this subsection we give the proof of Lemma \ref{lemma
4.1.2}.

\proof[Proof of Lemma \ref{lemma 4.1.2}]

First we claim that
\begin{equation}\label{equ 4.1.1}
\lim_{d\rightarrow+\infty}P_{\lambda,d}\Big(\sum_{m=0}^{\lfloor\sigma\log
d\rfloor}|V_m|>d^{0.2}\Big)=0.
\end{equation}
Equation \eqref{equ 4.1.1} follows from the following analysis. For
a given oriented path $\vec{l}:O=x_0\rightarrow
x_1\rightarrow\ldots\rightarrow x_m$ on $\mathbb{Z}^d_+$,
\[
P_{\lambda,d}\Big(\widetilde{U}(x_j,x_{j+1})<\widetilde{Y}(x_j)\text{~for
all~}0\leq j\leq m-1\Big)\leq \Big(\frac{\lambda}{d}M^2\Big)^m,
\]
since $\widetilde{Y}(\cdot)$ is an exponential time with rate one
while $\widetilde{U}(\cdot,\cdot)$ is an exponential time with rate
at most $\frac{\lambda}{d}M^2$. The number of oriented paths
starting at $O$ with length $m$ on $\mathbb{Z}^+_d$ is $d^m$. As a
result,
\[
E_{\lambda,d}|V_m|\leq
\Big(\frac{\lambda}{d}M^2\Big)^md^m=\big(\lambda M^2\big)^m,
\]
since $x\in V_m$ when and only when there exists an oriented path
$\vec{l}:O=x_0\rightarrow x_1\rightarrow\ldots\rightarrow x_m=x$
that $\widetilde{U}(x_j,x_{j+1})<\widetilde{Y}(x_j)$ for all $0\leq
j\leq m-1$.

Then, according to the Chebyshev's inequality and the fact that
$\sigma \log (\lambda M^2)<\frac{1}{10}$,
\[
P_{\lambda,d}\Big(\sum_{m=0}^{\lfloor\sigma\log
d\rfloor}|V_m|>d^{0.2}\Big)\leq
d^{-0.2}\sum_{m=0}^{\lfloor\sigma\log d\rfloor}\big(\lambda
M^2\big)^m\leq \frac{d^{-0.1}}{\lambda M^2-1},
\]
Equation \eqref{equ 4.1.1} follows from which directly.

For $d\geq 2$, we denote by $J(d)$ the event that
$\widetilde{U}(x,y)>\widetilde{Y}(x)$ for any $x\in
\bigcup_{m=0}^{\lfloor\sigma\log d\rfloor}V_m$ and any $y\in q(x)$.
It is easy to check that there exists a vertex $y$ satisfying
$x\rightarrow y,z\rightarrow y$ for given $x,z\in V_m$ when and only
when $x-z=e_i-e_j$ for some $1\leq i,j\leq d$ and such $y$ is unique
that $y=x+e_j=z+e_i$. Hence,
\begin{equation}\label{equ 4.1.2}
|q(x)|\leq |V_m|-1<|V_m|
\end{equation}
for any $x\in V_m$. By Equation \eqref{equ 4.1.2} and the fact that
$\widetilde{Y}(\cdot)$ is an exponential time with rate one while
$\widetilde{U}(\cdot,\cdot)$ is an exponential time with rate at
most $\frac{\lambda}{d}M^2$ while the event
$\sum_{m=0}^{\lfloor\sigma\log d\rfloor}|V_m|\leq d^{0.2}$ is
negative correlated with the event
$\widetilde{U}(x,y)<\widetilde{Y}(x)$ for $x\in
\bigcup_{m=0}^{\lfloor\sigma\log d\rfloor}V_m, y\in q(x)$, we have
\begin{align}\label{equ 4.1.3}
&P_{\lambda,d}\Big(J(d)^c\Bigg|\sum_{m=0}^{\lfloor\sigma\log
d\rfloor}|V_m|\leq d^{0.2}\Big) \notag\\
&\leq \sum_{m=0}^{\lfloor\sigma\log
d\rfloor}E_{\lambda,d}\Big(|V_m|\times|V_m|\frac{\lambda}{d}M^2\Bigg|\sum_{m=0}^{\lfloor\sigma\log
d\rfloor}|V_m|\leq d^{0.2}\Big)\\
&\leq \frac{\lambda
M^2}{d}E_{\lambda,d}\Big(\big(\sum_{m=0}^{\lfloor\sigma\log
d\rfloor}|V_m|\big)^2\Bigg|\sum_{m=0}^{\lfloor\sigma\log
d\rfloor}|V_m|\leq d^{0.2}\Big)\notag\\
&\leq \frac{\lambda M^2}{d}d^{0.4}=\lambda M^2 d^{-0.6}.\notag
\end{align}
By Equations \eqref{equ 4.1.1} and \eqref{equ 4.1.3},
\begin{equation}\label{equ 4.1.4}
\lim_{d\rightarrow+\infty}P_{\lambda,d}\Big(J(d),\sum_{m=0}^{\lfloor\sigma\log
d\rfloor}|V_m|\leq d^{0.2}\Big)=1.
\end{equation}
For $k<\lfloor \sigma\log d\rfloor$, conditioned on the event
\[
J(d)\bigcap \big\{\sum_{m=0}^{\lfloor\sigma\log d\rfloor}|V_m|\leq
d^{0.2}\big\}\bigcap\big\{\text{~the coupling is successful at
step~}k\big\},
\]
the coupling will be successful at step $k+1$ if
$U(g_k(x),y)>Y\big(g_k(x)\big)$ for any $x\in V_k$ and any $y$ that
$g_k(x)\Rightarrow y$ while $y\neq w_1,w_2,\ldots,w_{h(x)}$. Then,
according to a similar analysis with which leads to Equation
\eqref{equ 4.1.3},
\begin{align}\label{equ 4.1.5}
&P_{\lambda,d}\Big(\text{~the coupling is successful at
step~}k+1\\
&\text{\quad}\Bigg|J(d),\sum_{m=0}^{\lfloor\sigma\log
d\rfloor}|V_m|\leq
d^{0.2},\text{~the coupling is successful at step~}k\Big)\notag\\
&\geq 1-\frac{\lambda
M^2}{d}E_{\lambda,d}\Big(|V_m|\times|V_m|\Bigg|\sum_{m=0}^{\lfloor\sigma\log
d\rfloor}|V_m|\leq
d^{0.2}\Big) \notag\\
&\geq 1-\lambda M^2d^{-0.6}.\notag
\end{align}
By Equation \eqref{equ 4.1.5} and the conditional probability
formula,
\begin{equation}\label{equ 4.1.6}
P_{\lambda,d}\Big(B(d)\Bigg|J(d),\sum_{m=0}^{\lfloor\sigma\log
d\rfloor}|V_m|\leq d^{0.2}\Big)\geq \Big(1-\lambda
M^2d^{-0.6}\Big)^{\lfloor\sigma\log d\rfloor}.
\end{equation}
Lemma \ref{lemma 4.1.2} follows from Equations \eqref{equ 4.1.4} and
\eqref{equ 4.1.6} directly since
\[
\lim_{d\rightarrow+\infty}\Big(1-\lambda
M^2d^{-0.6}\Big)^{\lfloor\sigma\log d\rfloor}=1.
\]

\qed

\subsection{Proof of Lemma \ref{lemma 4.1}}\label{subsection 4.2}
In this subsection we give the proof of Lemma \ref{lemma 4.1}.

\proof[Proof of Lemma \ref{lemma 4.1}]

There is an isomorphism $\Phi:\mathbb{Z}_+^{d-\lfloor
\frac{d}{N(d)}\rfloor}\rightarrow \Gamma_1$ that
\[
\Phi(x_1,x_2,\ldots,x_{d-\lfloor
\frac{d}{N(d)}\rfloor})=(x_1,x_2,\ldots,x_{d-\lfloor
\frac{d}{N(d)}\rfloor},0,\ldots,0)
\]
for each $x=(x_1,x_2,\ldots,x_{d-\lfloor \frac{d}{N(d)}\rfloor})\in
\mathbb{Z}_+^{d-\lfloor \frac{d}{N(d)}\rfloor}$, hence
$\{\widehat{V}_n\}_{n\geq 0}$ can be identified as the SIR model
$\{V_n\}_{n\geq 0}$ on
$\mathbb{Z}_+^{d-\lfloor\frac{d}{N(d)}\rfloor}$ with infection rate
$\frac{\lambda(d-\lfloor \frac{d}{N(d)}\rfloor)}{d}$. For given
$\widetilde{\lambda}\in (\frac{1}{E(\rho^2)},\lambda)$ and
$\widetilde{\sigma}\in
(\sigma,\frac{1}{10\log(M^2\widetilde{\lambda})})$,
\[
\frac{\lambda(d-\lfloor \frac{d}{N(d)}\rfloor)}{d}\geq
\widetilde{\lambda} \text{~and~}\sigma \log d\leq
\widetilde{\sigma}\log (d-\lfloor\frac{d}{N(d)}\rfloor)
\]
for sufficiently large $d$. Then,
\[
P_{\lambda,d}\Big(\widehat{V}_{\lfloor\sigma\log d\rfloor}\neq
\emptyset\Big)\geq
P_{_{\widetilde{\lambda},d-\lfloor\frac{d}{N(d)}\rfloor}}\Big(V_{\lfloor\widetilde{\sigma}\log
(d-\lfloor\frac{d}{N(d)}\rfloor)\rfloor}\neq \emptyset\Big)
\]
for sufficiently large $d$ and hence
\[
\liminf_{d\rightarrow+\infty}P_{\lambda,d}\Big(\widehat{V}_{\lfloor\sigma\log
d\rfloor}\neq \emptyset\Big)\geq
\liminf_{d\rightarrow+\infty}P_{_{\widetilde{\lambda},d-\lfloor\frac{d}{N(d)}\rfloor}}\Big(V_{\lfloor\widetilde{\sigma}\log
(d-\lfloor\frac{d}{N(d)}\rfloor)\rfloor}\neq \emptyset\Big)\geq
E\big(\frac{\widetilde{\lambda}\rho\widetilde{\theta}}{1+\widetilde{\lambda}\rho\widetilde{\theta}}\big)
\]
according to Lemma \ref{lemma 4.1.1}, where $\widetilde{\theta}$
satisfies
\[
E\big(\frac{\widetilde{\lambda}\rho^2}{1+\widetilde{\lambda}\rho\widetilde{\theta}}\big)=1
\]
and it is easy to check that
$\lim_{\widetilde{\lambda}\rightarrow\lambda}\widetilde{\theta}=\theta$.

Let $\widetilde{\lambda}\rightarrow \lambda$, then the proof is
complete.

\qed

\subsection{Proof of Lemma \ref{lemma 4.2}}\label{subsection 4.3}
In this section we give the proof of Lemma \ref{lemma 4.2}. First we
introduce some notations and definitions.

We let
$\widehat{Y}_1,\widehat{Y}_2,\ldots,\widehat{Y}_{\lfloor\sigma \log
d\rfloor}$ be exponential times with rate one while
$\Lambda_1,\ldots,\Lambda_{\lfloor\sigma \log d\rfloor}$ be
exponential times with rate $\lambda M^2$. For $1\leq i\leq
\lfloor\sigma\log d\rfloor$ and $1\leq j\leq
\lfloor\frac{d}{N(d)}\rfloor$, let $\rho_{ij}$ be an independent
copy of $\rho$ while $\widehat{U}_{ij}$ is an exponential time with
rate $\frac{\lambda}{d}\epsilon\rho_{ij}$, where $\epsilon$ is
defined as in Equation \eqref{equ 4.1}. According to the basic
technique of measure theory, we can assume that $\{\rho_{ij}:1\leq
i\leq \lfloor\sigma\log d\rfloor, 1\leq j\leq
\lfloor\frac{d}{N(d)}\rfloor\}$ and $\{\rho(x)\}_{x\in
\mathbb{Z}^d_+}$ are defined under the same space and independent
under the annealed measure $P_{\lambda,d}$ while assume that
$\widetilde{U}(\cdot,\cdot),\widehat{U}_{\cdot\cdot},\widetilde{Y}(\cdot),\widehat{Y}_{\cdot},\Lambda_{\cdot}$
are defined under the same space and independent under the quenched
measure $P_{\lambda,\omega}$.

\begin{lemma}\label{lemma 4.3.1}
Let
$\xi_i=\sum\limits_{j=1}^{\lfloor\frac{d}{N(d)}\rfloor}1_{\{\widehat{U}_{ij}<Y_i\}}$
for $1\leq i\leq \lfloor\sigma\log d\rfloor$, then
\begin{align}\label{equ 4.3.1}
&P_{\lambda,d}\Big(|D|\geq K(d)\Big|\widehat{V}_{\lfloor\sigma\log
d\rfloor}\neq \emptyset\Big)\\
&\geq P_{\lambda,d}\Big(\sum_{i=1}^{\lfloor\sigma\log
d\rfloor}\xi_i\geq K(d)\Big|\widehat{Y}_i<\Lambda_i\text{~for
all~}i\leq \lfloor \sigma\log d\rfloor\Big), \notag
\end{align}
and
\begin{equation}\label{equ 4.3.2}
\lim_{d\rightarrow+\infty}E_{\lambda,d}\Big(e^{-s\frac{N(d)}{\log
d}\sum\limits_{j=1}^{\lfloor\sigma\log
d\rfloor}\xi_j}\Big|\widehat{Y}_i<\Lambda_i\text{~for all~}i\leq
\lfloor \sigma\log d\rfloor\Big)\exists=\Theta(s)
\end{equation}
for any $s>0$, where $\lim_{s\rightarrow+\infty}\Theta(s)=0$.
\end{lemma}

We give the proof of Lemma \ref{lemma 4.3.1} at the end of this
subsection. Now we show that how to utilize Lemma \ref{lemma 4.3.1}
to prove Lemma \ref{lemma 4.2}.

\proof[Proof of Lemma \ref{lemma 4.2}]

By Chebyshev's inequality, for any $s>0$,
\begin{align*}
&P_{\lambda,d}\Big(\sum_{i=1}^{\lfloor\sigma\log d\rfloor}\xi_i\leq
K(d)\Big|\widehat{Y}_i<\Lambda_i\text{~for all~}i\leq \lfloor
\sigma\log d\rfloor\Big)\\
&\leq e^{s\frac{K(d)N(d)}{\log
d}}E_{\lambda,d}\Big(e^{-s\frac{N(d)}{\log
d}\sum\limits_{j=1}^{\lfloor\sigma\log
d\rfloor}\xi_j}\Big|\widehat{Y}_i<\Lambda_i\text{~for all~}i\leq
\lfloor \sigma\log d\rfloor\Big).
\end{align*}
Then, according to Equation \eqref{equ 4.3.2},
\[
\limsup_{d\rightarrow+\infty}P_{\lambda,d}\Big(\sum_{i=1}^{\lfloor\sigma\log
d\rfloor}\xi_i\leq K(d)\Big|\widehat{Y}_i<\Lambda_i\text{~for
all~}i\leq \lfloor \sigma\log d\rfloor\Big)\leq \Theta(s)
\]
for any $s>0$, since $\lim_{d\rightarrow+\infty}\frac{K(d)N(d)}{\log
d}=0$.

Let $s\rightarrow+\infty$, since
$\lim_{s\rightarrow+\infty}\Theta(s)=0$, we have
\begin{equation}\label{equ 4.3.3}
\lim_{d\rightarrow+\infty}P_{\lambda,d}\Big(\sum_{i=1}^{\lfloor\sigma\log
d\rfloor}\xi_i\leq K(d)\Big|\widehat{Y}_i<\Lambda_i\text{~for
all~}i\leq \lfloor \sigma\log d\rfloor\Big)=0.
\end{equation}
Lemma \ref{lemma 4.2} follows from Equations \eqref{equ 4.3.1} and
\eqref{equ 4.3.3} directly.

\qed

Now we give the proof of Lemma \ref{lemma 4.3.1}.

\proof[Proof of Lemma \ref{lemma 4.3.1}]

Equation \eqref{equ 4.3.1} follows from the following analysis.
Conditioned on $\widehat{V}_{\lfloor\sigma\log d\rfloor}\neq
\emptyset$, there are $O\rightarrow X_1\rightarrow\ldots\rightarrow
X_{\lfloor\sigma\log d\rfloor}$ that $X_i\in \widehat{V}_i$ for
$1\leq i\leq \lfloor\sigma\log d\rfloor$. We choose
$X_1,\ldots,X_{\lfloor\sigma\log d\rfloor}$ as follows. We let
$X_{\lfloor\sigma\log d\rfloor}$ be the smallest one of
$\widehat{V}_{\lfloor\sigma\log d\rfloor}$ under the partial
$\prec$. Then, we let $X_{\lfloor\sigma\log d\rfloor-1}$ be the
smallest one of $\{x\in \widehat{V}_{\lfloor\sigma\log
d\rfloor-1}:x\rightarrow X_{\lfloor\sigma\log
d\rfloor}\text{~and~}\widetilde{U}(x,X_{\lfloor\sigma\log
d\rfloor})<\widetilde{Y}(x)\}$. By induction, if $X_j$ is well
defined for $i+1\leq j\leq \lfloor\sigma\log d\rfloor$, then we
define $X_i$ as the smallest one of
\[
\{x\in V_i:x\rightarrow
X_{i+1}\text{~and~}\widetilde{U}(x,X_{i+1})<\widetilde{Y}(x)\}.
\]

For $1\leq i\leq \lfloor\sigma\log d\rfloor$, let
\[
\zeta_i=\sum_{j=1}^{\lfloor\frac{d}{N(d)}\rfloor}1_{\{\widetilde{U}(X_i,X_i+e_{j+d-\lfloor\frac{d}{N(d)}\rfloor})<\widetilde{Y}(X_i)\}}
\text{~and~}
\eta_i=\sum_{j=1}^{\lfloor\frac{d}{N(d)}\rfloor}1_{\{\widehat{U}_{ij}<\widetilde{Y}(X_i)\}},
\]
then
\begin{equation}\label{equ 4.3.4}
|D|\geq \sum_{i=1}^{\lfloor\sigma\log d\rfloor}\zeta_i
\end{equation}
according to our definition of $D$.

Since $X_i\in \widehat{V}_i$, $\rho(X_i)>0$ and hence $\rho(X_i)\geq
\epsilon$ by assumption \eqref{equ 4.1}. As a result,
$\widetilde{U}(X_i,X_i+e_{j+d-\lfloor\frac{d}{N(d)}\rfloor})$ is an
exponential time with rate at least
$\frac{\lambda}{d}\epsilon\rho(x_i+e_{j+d-\lfloor\frac{d}{N(d)}\rfloor})$,
where $\rho(X_i+e_{j+d-\lfloor\frac{d}{N(d)}\rfloor})$ is an
independent copy of $\rho$. Therefore,
$\widetilde{U}(X_i,X_i+e_{j+d-\lfloor\frac{d}{N(d)}\rfloor})$ is
stochastic dominated from above by $\widehat{U}_{ij}$ and $\eta_i$
is dominated from above by $\xi_i$. Hence,
\begin{equation}\label{equ 4.3.5}
P_{\lambda,d}\Big(\sum_{i=1}^{\lfloor\sigma\log d\rfloor}\zeta_i\geq
K(d)\Big|\widehat{V}_{\lfloor\sigma\log d\rfloor}\neq \emptyset\Big)
\geq P_{\lambda,d}\Big(\sum_{i=1}^{\lfloor\sigma\log
d\rfloor}\eta_i\geq K(d)\Big|\widehat{V}_{\lfloor\sigma\log
d\rfloor}\neq \emptyset\Big).
\end{equation}
For any oriented path $\vec{l}:O\rightarrow
l_1\rightarrow\ldots\rightarrow l_{\lfloor\sigma\log d\rfloor}$, we
use $\gamma(\vec{l})$ to denote
\[
P_{\lambda,d}\Big(X_i=l_i\text{~for all~}1\leq i\leq
\lfloor\sigma\log d\rfloor\Bigg|\widehat{V}_{\lfloor\sigma\log
d\rfloor}\neq \emptyset\Big).
\]
Then,
\begin{align}\label{equ 4.3.6}
&P_{\lambda,d}\Big(\sum_{i=1}^{\lfloor\sigma\log d\rfloor}\eta_i\geq
K(d)\Big|\widehat{V}_{\lfloor\sigma\log d\rfloor}\neq
\emptyset\Big)\\
&=\sum_{\vec{l}}\gamma(\vec{l})P_{\lambda,d}\Big(\sum_{i=1}^{\lfloor\sigma\log
d\rfloor}\eta_i(\vec{l})\geq K(d)\Big|X_i=l_i\text{~for all~}1\leq
i\leq \lfloor\sigma\log d\rfloor, \widehat{V}_{\lfloor\sigma\log
d\rfloor}\neq \emptyset\Big), \notag
\end{align}
where
\[
\eta_i(\vec{l})=\sum_{j=1}^{\lfloor\frac{d}{N(d)}\rfloor}1_{\{\widehat{U}_{ij}<\widetilde{Y}(l_i)\}}.
\]
The condition $\{X_i=l_i\text{~for all~}1\leq i\leq
\lfloor\sigma\log d\rfloor, \widehat{V}_{\lfloor\sigma\log
d\rfloor}\neq \emptyset\}$ in Equation \eqref{equ 4.3.6} is
concerned with the values of $\widetilde{Y}(l_i)$ and
$\{\widetilde{U}(l_i,l_i+e_j):1\leq j\leq
d-\lfloor\frac{d}{N(d)}\rfloor\}$ for $1\leq i\leq \lfloor\sigma\log
d\rfloor$. A worse condition for $\sum_{i=1}^{\lfloor\sigma\log
d\rfloor}\eta_i(\vec{l})\geq K(d)$ to occur is that
\[
\widetilde{Y}(l_i)<\inf\big\{\widetilde{U}(l_i,l_i+e_j):1\leq j\leq
d-\lfloor\frac{d}{N(d)}\rfloor\big\}
\]
for all $1\leq i\leq \lfloor\sigma\log d\rfloor$, i.e.,
\begin{align}\label{equ 4.3.7}
&P_{\lambda,d}\Big(\sum_{i=1}^{\lfloor\sigma\log
d\rfloor}\eta_i(\vec{l})\geq K(d)\Big|X_i=l_i\text{~for all~}1\leq
i\leq \lfloor\sigma\log d\rfloor, \widehat{V}_{\lfloor\sigma\log
d\rfloor}\neq \emptyset\Big)\geq \notag\\
&P_{\lambda,d}\Big(\sum_{i=1}^{\lfloor\sigma\log
d\rfloor}\eta_i(\vec{l})\geq
K(d)\Big|\widetilde{Y}(l_i)<\inf\big\{\widetilde{U}(l_i,l_i+e_j):1\leq
j\leq d-\lfloor\frac{d}{N(d)}\rfloor\big\} \notag\\
&\text{~for all~}1\leq i\leq \lfloor\sigma\log d\rfloor\Big).
\end{align}
Note that $\widetilde{Y}(l_i)$ is with the same distribution as that
of $\widehat{Y}_i$ while $\eta_i(\vec{l})$ is with the same
distribution as that of $\xi_i$. Further more,
$\inf\big\{\widetilde{U}(l_i,l_i+e_j):1\leq j\leq
d-\lfloor\frac{d}{N(d)}\rfloor\big\}$ is an exponential time with
rate at most
\[
(d-\lfloor\frac{d}{N(d)}\rfloor)\frac{\lambda}{d} M^2\leq \lambda
M^2,
\]
which is the rate of $\Lambda_i$. Then,
\begin{align*}
&P_{\lambda,d}\Big(\sum_{i=1}^{\lfloor\sigma\log
d\rfloor}\eta_i(\vec{l})\geq
K(d)\Big|\widetilde{Y}(l_i)<\inf\big\{\widetilde{U}(l_i,l_i+e_j):1\leq
j\leq d-\lfloor\frac{d}{N(d)}\rfloor\big\} \notag\\
&\text{~for all~}1\leq i\leq \lfloor\sigma\log d\rfloor\Big)\\
&\geq P_{\lambda,d}\Big(\sum_{i=1}^{\lfloor\sigma\log
d\rfloor}\xi_i\geq K(d)\Big|\widehat{Y}_i<\Lambda_i\text{~for
all~}i\leq \lfloor \sigma\log d\rfloor\Big).
\end{align*}
Therefore, by Equation \eqref{equ 4.3.7},
\begin{align}\label{equ 4.3.8}
&P_{\lambda,d}\Big(\sum_{i=1}^{\lfloor\sigma\log
d\rfloor}\eta_i(\vec{l})\geq K(d)\Big|X_i=l_i\text{~for all~}1\leq
i\leq \lfloor\sigma\log d\rfloor, \widehat{V}_{\lfloor\sigma\log
d\rfloor}\neq \emptyset\Big)\notag\\
&\geq P_{\lambda,d}\Big(\sum_{i=1}^{\lfloor\sigma\log
d\rfloor}\xi_i\geq K(d)\Big|\widehat{Y}_i<\Lambda_i\text{~for
all~}i\leq \lfloor \sigma\log d\rfloor\Big).
\end{align}

By Equations \eqref{equ 4.3.6} and \eqref{equ 4.3.8},
\begin{align}\label{equ 4.3.9}
&P_{\lambda,d}\Big(\sum_{i=1}^{\lfloor\sigma\log d\rfloor}\eta_i\geq
K(d)\Big|\widehat{V}_{\lfloor\sigma\log d\rfloor}\neq
\emptyset\Big)\\
&\geq \sum_{\vec{l}}\gamma(\vec{l})
P_{\lambda,d}\Big(\sum_{i=1}^{\lfloor\sigma\log d\rfloor}\xi_i\geq
K(d)\Big|\widehat{Y}_i<\Lambda_i\text{~for all~}i\leq \lfloor
\sigma\log d\rfloor\Big)\notag\\
&=P_{\lambda,d}\Big(\sum_{i=1}^{\lfloor\sigma\log d\rfloor}\xi_i\geq
K(d)\Big|\widehat{Y}_i<\Lambda_i\text{~for all~}i\leq \lfloor
\sigma\log d\rfloor\Big),\notag
\end{align}
since $\sum_{\vec{l}}\gamma(\vec{l})=1$. Equation \eqref{equ 4.3.1}
follows directly from Equations \eqref{equ 4.3.4}, \eqref{equ 4.3.5}
and \eqref{equ 4.3.9}.

Equation \eqref{equ 4.3.2} follows from the following analysis.
According to the assumption of independence of the exponential
times,
\begin{align}\label{equ 4.3.10}
&E_{\lambda,d}\Big(e^{-s\frac{N(d)}{\log
d}\sum\limits_{j=1}^{\lfloor\sigma\log
d\rfloor}\xi_j}\Big|\widehat{Y}_i<\Lambda_i\text{~for all~}i\leq
\lfloor \sigma\log d\rfloor\Big)\notag\\
&=\Bigg[E_{\lambda,d}\Big(e^{-s\frac{N(d)}{\log
d}\xi_1}\Bigg|\widehat{Y}_1<\Lambda_1\Big)\Bigg]^{\lfloor\sigma\log
d\rfloor}.
\end{align}
Let $Y_0$ be an exponential time with rate one while $\Lambda_0$ be
an exponential time with rate $\lambda M^2$ and independent of
$Y_0$, then by direct calculation,
\[
E_{\lambda,d}\Big(e^{-s\frac{N(d)}{\log
d}\xi_1}\Bigg|\widehat{Y}_1<\Lambda_1\Big)=E\Big(\Xi(\widehat{Y}_0,s,d)\Big|\widehat{Y}_0<\Lambda_0\Big),
\]
where
\[
\Xi(t,s,d)=\Bigg[E\Big(e^{-s\frac{N(d)}{\log
d}}\big(1-e^{-\frac{\lambda t\epsilon\rho
}{d}}\big)+e^{-\frac{\lambda
t\epsilon\rho}{d}}\Big)\Bigg]^{\lfloor\frac{d}{N(d)}\rfloor}.
\]
Then, by Equation \eqref{equ 4.3.10},
\begin{align}\label{equ 4.3.11}
&E_{\lambda,d}\Big(e^{-s\frac{N(d)}{\log
d}\sum\limits_{j=1}^{\lfloor\sigma\log
d\rfloor}\xi_j}\Big|\widehat{Y}_i<\Lambda_i\text{~for all~}i\leq
\lfloor \sigma\log d\rfloor\Big)\notag\\
&=\Bigg[E\Big(\Xi(\widehat{Y}_0,s,d)\Big|\widehat{Y}_0<\Lambda_0\Big)\Bigg]^{\lfloor\sigma\log
d\rfloor}.
\end{align}
By direct calculation, it is not difficult to check that
\[
\lim_{d\rightarrow+\infty}\lfloor\sigma\log
d\rfloor\Big(\Xi(t,s,d)-1\Big)=-s\sigma\lambda t\epsilon E\rho.
\]
for any $t>0$. As a result,
\begin{equation}\label{equ 4.3.12}
\lim_{d\rightarrow+\infty}\Bigg[E\Big(\Xi(\widehat{Y}_0,s,d)\Big|\widehat{Y}_0<\Lambda_0\Big)\Bigg]^{\lfloor\sigma\log
d\rfloor}=E\Big(e^{-s\sigma\lambda \widehat{Y}_0\epsilon
E\rho}\Big|\widehat{Y}_0<\Lambda_0\Big):=\Theta(s).
\end{equation}
Note that here we still use the theorem of calculus that
$(1+a_d)^{c_d}\rightarrow e^c$ when $a_d\rightarrow
0,c_d\rightarrow+\infty$ and $a_dc_d\rightarrow c$. Equation
\eqref{equ 4.3.2} follows directly from Equations \eqref{equ 4.3.11}
and \eqref{equ 4.3.12}.

\qed

\subsection{Proof of Lemma \ref{lemma 4.3}}\label{subsection 4.4}
In this subsection we give the proof of Lemma \ref{lemma 4.3}. The
proof is inspired a lot by the approach introduced in
\cite{Xue2015}. First we introduce some definitions and notations.
We let $\{\vartheta_n\}_{n\geq 0}$ be the oriented random walk on
$\mathbb{Z}_+^d$ such that
\[
P(\vartheta_{n+1}-\vartheta_n=e_i)=\frac{1}{d}
\]
for each $n\geq 0$ and $1\leq i\leq d$. We let $\{\nu_n\}_{n\geq 0}$
be an independent copy of $\{\vartheta_n\}_{n\geq 0}$. From now on,
we denote by $\mathbb{P}$ the probability measure of
$\{\vartheta_n\}_{n\geq 0}$ and $\{\nu_n\}_{n\geq 0}$ while denote
by $\mathbb{E}$ the expectation with respect to $\mathbb{P}$. When
we need to point out the dimension $d$ of the lattice, we write
$\mathbb{P}$ and $\mathbb{E}$ as $\mathbb{P}_d$ and $\mathbb{E}_d$.
We write $\vartheta_n$ (resp. $\nu_n$) as $\vartheta_n^x$ (resp.
$\nu_n^x$) when $\vartheta_0=x$ (resp. $\nu_0=x$). For $x,y\in
\mathbb{Z}^d_+$ that $x\neq y$ and $\|x\|\leq \|y\|$, we define
\[
\tau_{x,y}=\inf\big\{k\geq
\|y\|-\|x\|:~\vartheta_k^x=\nu^y_{k-\|y\|+\|x\|}\big\}.
\]
That is to say, $\tau_{x,y}$ is the first moment when
$\{\vartheta_n^x\}_{n\geq 0}$ visits some vertex on the path of
$\{\nu_n^y\}_{n\geq 0}$. Note that
\[
\|\vartheta_n\|=\|\vartheta_0\|+n\text{~and~}\|\nu_n\|=\|\nu_0\|+n
\]
according to the definition of the oriented random walk, hence
\[
k=l+\|y\|-\|x\|
\]
when $\vartheta_k^x=\nu_l^y$ for some $k,l$. For $x,y\in
\mathbb{Z}_+^d$ that $\|x\|\leq \|y\|$, we introduce the following
random variables. We define
\[
\tau_0^{x,y}=
\begin{cases}
0  & \text{~if~} x=y,\\
\tau_{x,y} & \text{~if~}x\neq y.
\end{cases}
\]
We let
\begin{align*}
& \tau_1^{x,y}=\inf\big\{n\geq
\tau_0^{x,y}:~\vartheta_n^x=\nu^y_{n-\|y\|+\|x\|},\vartheta_{n+1}^x=\nu^y_{n+1-\|y\|+\|x\|}\big\},\\
& \kappa_1^{x,y}=\inf\big\{n>\tau_1^{x,y}:~\vartheta_n^x=\nu^y_{n-\|y\|+\|x\|},\vartheta_{n+1}^x\neq \nu^y_{n+1-\|y\|+\|x\|}\big\}, \\
& \tau_2^{x,y}=\inf\big\{n>
\kappa_1^{x,y}:~\vartheta_n^x=\nu^y_{n-\|y\|+\|x\|},\vartheta_{n+1}^x=\nu^y_{n+1-\|y\|+\|x\|}\big\},\\
& \kappa_2^{x,y}=\inf\big\{n>\tau_2^{x,y}:~\vartheta_n^x=\nu^y_{n-\|y\|+\|x\|},\vartheta_{n+1}^x\neq \nu^y_{n+1-\|y\|+\|x\|}\big\}, \\
&\ldots\ldots\\
& \tau_l^{x,y}=\inf\big\{n>
\kappa_{l-1}^{x,y}:~\vartheta_n^x=\nu^y_{n-\|y\|+\|x\|},\vartheta_{n+1}^x=\nu^y_{n+1-\|y\|+\|x\|}\big\},\\
&
\kappa_l^{x,y}=\inf\big\{n>\tau_l^{x,y}:~\vartheta_n^x=\nu^y_{n-\|y\|+\|x\|},\vartheta_{n+1}^x\neq
\nu^y_{n+1-\|y\|+\|x\|}\big\},\\
&\ldots\ldots
\end{align*}
That is to say, $\tau_1^{x,y}$ is the first moment $n$ such that
$\vartheta_n^x=\nu_{n-\|y\|+\|x\|}^y$ and
$\vartheta_{n+1}^x=\nu_{n+1-\|y\|+\|x\|}^y$. For $l\geq 1$,
$\kappa_l^{x,y}$ is the first moment $n$ after $\tau_l^{x,y}$ such
that $\vartheta_n^x=\nu_{n-\|y\|+\|x\|}^y$ and
$\vartheta_{n+1}^x\neq \nu_{n+1-\|y\|+\|x\|}^y$ while
$\tau_{l+1}^{x,y}$ is the first moment $n$ after $\kappa_l^{x,y}$
that $\vartheta_n^x=\nu_{n-\|y\|+\|x\|}^y$ and
$\vartheta_{n+1}^x=\nu_{n+1-\|y\|+\|x\|}^y$.

We define
\[
T(x,y)=
\begin{cases}
\sup\big\{l\geq 0:\tau_l^{x,y}<+\infty\big\} & \text{~if~}
\tau_0^{x,y}<+\infty,\\
0 & \text{~if~}\tau_0^{x,y}=+\infty.
\end{cases}
\]
In this subsection we assume that $d\geq 4$ such that
$T(x,y)<+\infty$ with probability one according to the conclusion
given in \cite{Cox1983} about the collision times of two independent
oriented random walks.

For $1\leq l\leq T(x,y)$, we define
\[
h_l^{x,y}=\kappa_l^{x,y}-\tau_l^{x,y}.
\]
We let
\begin{align*}
&f_0^{x,y}=\big|\{\tau_0^{x,y}\leq n<\tau_1^{x,y}:\vartheta_n^x=\nu_{n-\|y\|+\|x\|}^y\}\big|,\\
&f_1^{x,y}=\big|\{\kappa_1^{x,y}<
n<\tau_2^{x,y}:\vartheta_n^x=\nu_{n-\|y\|+\|x\|}^y\}\big|,\\
&\ldots\ldots\\
&f_l^{x,y}=\big|\{\kappa_l^{x,y}<
n<\tau_{l+1}^{x,y}:\vartheta_n^x=\nu_{n-\|y\|+\|x\|}^y\}\big|,\\
&\ldots\ldots\\
&f_{T(x,y)-1}^{x,y}=\big|\{\kappa_{T(x,y)-1}^{x,y}<
n<\tau_{T(x,y)}^{x,y}:\vartheta_n^x=\nu_{n-\|y\|+\|x\|}^y\}\big|,\\
&f_{T(x,y)}^{x,y}=\big|\{n>\kappa_{T(x,y)}^{x,y}
:\vartheta_n^x=\nu_{n-\|y\|+\|x\|}^y\}\big|,
\end{align*}
where $|A|$ is the cardinality of the set $A$ as we have introduced.
Then, for $x, y$ that $\|x\|\leq \|y\|$, we define
\begin{align*}
&R(x,y)=\\
&\frac{2^{T(x,y)+\sum\limits_{i=0}^{T(x,y)}f_i^{x,y}}\big(1+\frac{\lambda
M^2}{d}\big)^{4T(x,y)+2\sum\limits_{i=1}^{T(x,y)}h_i^{x,y}+4\sum\limits_{i=0}^{T(x,y)}f_i^{x,y}}M^{6T(x,y)+4\sum\limits_{i=0}^{T(x,y)}f_i^{x,y}}}
{\Big(\frac{\lambda
E(\rho^2)}{d}\Big)^{\sum\limits_{i=1}^{T(x,y)}h_i^{x,y}}\Big(E(\rho^2)\Big)^{3T(x,y)+2\sum\limits_{i=0}^{T(x,y)}f_i(x,y)}}
\end{align*}
when $\tau_0^{x,y}>\|y\|-\|x\|$ while define
\begin{align*}
&R(x,y)=\\
&\frac{2^{T(x,y)+\sum\limits_{i=1}^{T(x,y)}f_i^{x,y}}\big(1+\frac{\lambda
M^2}{d}\big)^{4T(x,y)+2\sum\limits_{i=1}^{T(x,y)}h_i^{x,y}+4\sum\limits_{i=1}^{T(x,y)}f_i^{x,y}-1}M^{6T(x,y)+4\sum\limits_{i=1}^{T(x,y)}f_i^{x,y}-1}}
{ \Big(\frac{\lambda
E(\rho^2)}{d}\Big)^{\sum\limits_{i=1}^{T(x,y)}h_i^{x,y}}\Big(E(\rho^2)\Big)^{3T(x,y)+2\sum\limits_{i=1}^{T(x,y)}f_i(x,y)-1}E\rho}
\end{align*}
when $\|y\|-\|x\|=\tau_0^{x,y}=\tau_1^{x,y}$ and define
\begin{align*}
&R(x,y)=\\
&\frac{2^{T(x,y)+\sum\limits_{i=0}^{T(x,y)}f_i^{x,y}}\big(1+\frac{\lambda
M^2}{d}\big)^{4T(x,y)+2\sum\limits_{i=1}^{T(x,y)}h_i^{x,y}+4\sum\limits_{i=0}^{T(x,y)}f_i^{x,y}-1}M^{6T(x,y)+4\sum\limits_{i=0}^{T(x,y)}f_i^{x,y}-1}}
{\Big(\frac{\lambda
E(\rho^2)}{d}\Big)^{\sum\limits_{i=1}^{T(x,y)}h_i^{x,y}}\Big(E(\rho^2)\Big)^{3T(x,y)+2\sum\limits_{i=0}^{T(x,y)}f_i(x,y)-1}E\rho}
\end{align*}
when $\|y\|-\|x\|=\tau_0^{x,y}<\tau_1^{x,y}$. For
$x,y$ that $\|x\|>\|y\|$, we define
\[
\tau_i^{x,y},\kappa_i^{x,y},h_i^{x,y},T(x,y),f_i^{x,y},R(x,y) \text{~as~}
\tau_i^{y,x},\kappa_i^{x,y},h_i^{y,x},T(y,x),f_i^{y,x},R(y,x).
\]

The following three lemmas is crucial for us to prove Lemma
\ref{lemma 4.3}.

\begin{lemma}\label{lemma 4.4.1}
There exists $c_1>0$ which does not depend on $d$ such that
\[
\mathbb{P}_d\big(\tau_{x,y}<+\infty\big)\leq \frac{c_1}{\sqrt{d}}
\]
for any $d\geq 4$, $x,y\in \mathbb{Z}^+_d,x\neq y$.
\end{lemma}

\begin{lemma}\label{lemma 4.4.2}
For given $\lambda>\frac{1}{E(\rho^2)}$, there exists $d_0\geq 4$
and $c_2>0$ which does not depend on $d$ such that
\[
\mathbb{E}_d\Big(R(x,y)\Big|\tau(x,y)<+\infty\Big)\leq c_2
\]
for any $d\geq d_0$, $x,y\in \mathbb{Z}^+_d$.
\end{lemma}

\begin{lemma}\label{lemma 4.4.3}
For $A\subseteq \mathbb{Z}_d^+$,
\begin{align*}
P_{\lambda,d}\Big(I_t^A\neq \emptyset, \forall~t\geq
0\Big)\geq \frac{1}{\frac{1}{|A|^2}\sum\limits_{x\in A}\sum\limits_{y\in A}\mathbb{E}_d\big(R(x,y)\big)}.
\end{align*}
\end{lemma}

The proofs of Lemmas \ref{lemma 4.4.1}-\ref{lemma 4.4.3} will be
given later. Now we give the proof of Lemma \ref{lemma 4.3}.

\proof[Proof of Lemma \ref{lemma 4.3}]
For $x\neq y$, according to the definition of $R(x,y)$, $R(x,y)=1$ when $\tau_{x,y}=+\infty$. Therefore,
by Lemmas \ref{lemma 4.4.1} and \ref{lemma 4.4.2},
\begin{equation}\label{equ 4.4.1}
\mathbb{E}_d\big(R(x,y)\big)=\mathbb{P}_d\big(\tau_{x,y}=+\infty\big)+\mathbb{E}_d\big(R(x,y)1_{\{\tau_{x,y}<+\infty\}}\big)\leq
1+\frac{c_1c_2}{\sqrt{d}}
\end{equation}
for any $d\geq d_0$, $x,y\in \mathbb{Z}^d_+,x\neq y$. By Lemma
\ref{lemma 4.4.2} and Equation \eqref{equ 4.4.1},
\begin{equation}\label{equ 4.4.one half}
\sum\limits_{x\in A}\sum\limits_{y\in A}\mathbb{E}_d\big(R(x,y)\big)\leq |A|c_2+(|A|^2-|A|)(1+\frac{c_1c_2}{\sqrt{d}}).
\end{equation}
By Lemma \ref{lemma 4.4.3} and Equation \eqref{equ 4.4.one half},
\begin{equation}\label{equ 4.4.2}
P_{\lambda,d}\Big(I_t^A\neq \emptyset, \forall~t\geq
0\Big)\geq
\frac{|A|^2}{(|A|^2-|A|)(1+\frac{c_1c_2}{\sqrt{d}})+|A|c_2}
\end{equation}
for any $d\geq d_0$, $A\subseteq \mathbb{Z}^d_+$.

Let
$m(d)=\frac{|K(d)|^2}{(|K(d)|^2-|K(d)|)(1+\frac{c_1c_2}{\sqrt{d}})+|K(d)|c_2}$,
then $\lim_{d\rightarrow+\infty}m(d)=1$ while
\[
P_{\lambda,d}\Big(I_t^A\neq \emptyset, \forall~t\geq
0\Big)\geq m(d)
\]
for any $A\subseteq \mathbb{Z}_+^d,|A|=K(d)$ by Equation \eqref{equ
4.4.2} and the proof is complete.

\qed

Now we give the proof of Lemma \ref{lemma 4.4.1}.

\proof[Proof of Lemma \ref{lemma 4.4.1}]

Let
\[
\tau_{_{O,O}}=\inf\big\{n\geq 1:\vartheta_n^O=\nu_n^O\big\},
\]
then by  the conclusion given in \cite{Cox1983}, there exists $c_3>0$ which does not depend on $d$ that
\begin{equation}\label{equ 4.4.3}
\mathbb{P}_d\big(\tau_{_{O,O}}<+\infty\big)\leq \frac{1}{d}+\frac{c_3}{d^2}
\end{equation}
for all $d\geq 4$. Since $\mathbb{P}_d\big(\tau_{_{O,O}}=1\big)=\frac{1}{d}$, according to the spatial homogeneity of $\mathbb{Z}^d_+$,
\begin{equation}\label{equ 4.4.4}
(1-\frac{1}{d})\mathbb{P}_d\big(\tau_{e_i,e_j}<+\infty\big)=\mathbb{P}_d\big(2\leq \tau_{_{O,O}}<+\infty\big)\leq \frac{c_3}{d^2}
\end{equation}
for any $d\geq 4, 1\leq i<j\leq d$.
For $x,y$ that $x\neq y$ and $\|x\|=\|y\|$, $\|x-y\|$ is an even at least two. Let
\[
\widehat{\tau}_{x,y}=\inf\{n\geq 0:\|\vartheta_n^x-\nu_n^y\|=2\},
\]
then, according to the strong Markov property,
\begin{equation}\label{equ 4.4.5}
\mathbb{P}_d\big(\tau_{x,y}<+\infty\big)=\mathbb{P}_d\big(\widehat{\tau}_{x,y}<+\infty\big)\mathbb{P}_d\big(\tau_{e_1,e_2}<+\infty\big)\leq \frac{c_3}{d(d-1)},
\end{equation}
since $\|\vartheta_{n+1}^x-\nu_{n+1}^y\|-\|\vartheta_n^x-\nu_n^y\|\in\{0,-2,2\}$ for each $n$.

For $x,y$ that $\|x\|<\|y\|$,
\begin{align}\label{equ 4.4.6}
&\mathbb{P}_d\big(\tau_{x,y}<+\infty\big)\notag\\
&=\mathbb{P}_d\big(\vartheta_{\|y\|-\|x\|}^x=\|y\|\big)
+\mathbb{P}_d\big(\tau_{x,y}<+\infty\big|\vartheta^x_{\|y\|-\|x\|}\neq\|y\|\big)\mathbb{P}_d\big(\vartheta^x_{\|y\|-\|x\|}\neq\|y\|\big)\notag\\
&\leq \mathbb{P}_d\big(\vartheta_{\|y\|-\|x\|}^x=\|y\|\big)
+\mathbb{P}_d\big(\tau_{x,y}<+\infty\big|\vartheta^x_{\|y\|-\|x\|}\neq\|y\|\big).
\end{align}
Since $\|\vartheta^x_{\|y\|-\|x\|}\|=\|y\|$, by Equation \eqref{equ 4.4.5} and the strong Markov property,
\begin{equation}\label{equ 4.4.7}
\mathbb{P}_d\big(\tau_{x,y}<+\infty\big|\vartheta^x_{\|y\|-\|x\|}\neq\|y\|\big)\leq \frac{c_3}{d(d-1)}.
\end{equation}
By Equation \eqref{equ 4.4.3}, there exists $c_4>0$ which does not depend on $d$ that
\begin{align*}
&\mathbb{P}_d\big(\vartheta_{\|y\|-\|x\|}^x=\|y\|\big)^2=\mathbb{P}_d\big(\vartheta_{\|y\|-\|x\|}^x=\|y\|\big)\mathbb{P}_d\big(\nu_{\|y\|-\|x\|}^x=\|y\|\big)\\
&\leq \mathbb{P}_d\big(\vartheta_{\|y\|-\|x\|}^x=\nu_{\|y\|-\|x\|}^x\big)=\mathbb{P}_d\big(\vartheta_{\|y\|-\|x\|}^O=\nu_{\|y\|-\|x\|}^O\big)\\
&\leq \mathbb{P}_d\big(\tau_{_{O,O}}<+\infty\big)\leq \frac{c_4}{d}.
\end{align*}
As a result,
\begin{equation}\label{equ 4.4.8}
\mathbb{P}_d\big(\vartheta_{\|y\|-\|x\|}^x=\|y\|\big)\leq \frac{\sqrt{c_4}}{\sqrt{d}}.
\end{equation}
By Equations \eqref{equ 4.4.6},\eqref{equ 4.4.7} and \eqref{equ 4.4.8},
\begin{equation}\label{equ 4.4.9}
\mathbb{P}_d\big(\tau_{x,y}<+\infty\big)\leq \frac{\sqrt{c_4}}{\sqrt{d}}+\frac{c_3}{d(d-1)}
\end{equation}
for any $d\geq 4$, $x,y\in \mathbb{Z}^d_+$ that $\|x\|<\|y\|$.
Lemma \ref{lemma 4.4.1} follows directly from Equations \eqref{equ 4.4.5} and \eqref{equ 4.4.9}.

\qed

Now we give the proof of Lemma \ref{lemma 4.4.2}.

\proof[Proof of Lemma \ref{lemma 4.4.2}]

According to the definition of $R(\cdot,\cdot)$, for any $x,y\in \mathbb{Z}^+_d$,
\begin{align}\label{equ 4.4.10}
&R(x,y)\leq \\
&\frac{c_52^{T(x,y)+\sum\limits_{i=0}^{T(x,y)}f_i^{x,y}}\big(1+\frac{\lambda
M^2}{d}\big)^{4T(x,y)+2\sum\limits_{i=1}^{T(x,y)}h_i^{x,y}+4\sum\limits_{i=0}^{T(x,y)}f_i^{x,y}}M^{6T(x,y)+4\sum\limits_{i=0}^{T(x,y)}f_i^{x,y}}}
{\Big(\frac{\lambda
E(\rho^2)}{d}\Big)^{\sum\limits_{i=1}^{T(x,y)}h_i^{x,y}}\Big(E(\rho^2)\Big)^{3T(x,y)+2\sum\limits_{i=0}^{T(x,y)}f_i(x,y)}}, \notag
\end{align}
where $c_5>0$ is a constant which depends on $M,E(\rho^2), E\rho,\lambda$ and does not depend on $d$. According to the strong Markov property,
for any positive integers $T,\{h_i\}_{i=1}^T,\{f_i\}_{i=0}^T$,
\begin{align}\label{equ 4.4.11}
&\mathbb{P}_d\Big(T(x,y)=T,h_i^{x,y}=h_i,f_i^{x,y}=f_i\text{~for~}0\leq i\leq T\Big|\tau_{x,y}<+\infty\Big)\\
&\leq \mathbb{P}_d\Big(2\leq \tau_{_{O,O}}<+\infty\Big)^{f_0+\sum_{i=1}^{T-1}(f_i+1)+f_T}\mathbb{P}_d\Big(\tau_{_{O,O}}=1\Big)^{\sum_{i=1}^Th_i}, \notag
\end{align}
where $\tau_{_{O,O}}$ is defined as in the proof of Lemma \ref{lemma 4.4.1}. By Equations \eqref{equ 4.4.4} and \eqref{equ 4.4.11},
\begin{align}\label{equ 4.4.12}
&\mathbb{P}_d\Big(T(x,y)=T,h_i^{x,y}=h_i,f_i^{x,y}=f_i\text{~for~}0\leq i\leq T\Big|\tau_{x,y}<+\infty\Big)\\
&\leq \Big(\frac{c_3}{d^2}\Big)^{f_0-1+\sum_{i=1}^Tf_i+T}\big(\frac{1}{d}\big)^{\sum_{i=1}^Th_i}. \notag
\end{align}
By Equations \eqref{equ 4.4.10} and \eqref{equ 4.4.12},
\begin{align}\label{equ 4.4.13}
&\mathbb{E}_d\Big(R(x,y)\Big|\tau_{x,y}<+\infty\Big) \\
&\leq c_5\sum_{T=0}^{+\infty}\sum_{f_0=1}^{+\infty}\sum_{f_1=1}^{+\infty}\ldots\sum_{f_T=1}^{+\infty}\sum_{h_1=1}^{+\infty}\ldots
\sum_{h_T=1}^{+\infty} \big(\frac{c_3}{d^2}\big)^{\sum_{i=0}^Tf_i+T-1}\big(\frac{1}{d}\big)^{\sum_{i=1}^Th_i}\notag\\
&\times \frac{2^{T+\sum\limits_{i=0}^{T}f_i}\big(1+\frac{\lambda
M^2}{d}\big)^{4T+2\sum\limits_{i=1}^{T}h_i+4\sum\limits_{i=0}^{T}f_i}M^{6T+4\sum\limits_{i=0}^{T}f_i}}
{\Big(\frac{\lambda
E(\rho^2)}{d}\Big)^{\sum\limits_{i=1}^{T}h_i}\Big(E(\rho^2)\Big)^{3T+2\sum\limits_{i=0}^{T}f_i}}\notag\\
&=\widehat{c}_5\sum_{T=0}^{+\infty}\big(c_7(d)\big)^T\Big[\sum_{f_0=0}^{+\infty}\sum_{f_1=1}^{+\infty}\ldots\sum_{f_T=1}^{+\infty}\big(c_8(d)\big)^{\sum_{i=0}^{T}f_i}\Big]\notag\\
&\times \Big[\sum_{h_1=1}^{+\infty}\ldots\sum_{h_T=1}^{+\infty}\big(c_9(d)\big)^{\sum_{i=1}^Th_i}\Big],\notag
\end{align}
where $\widehat{c}_5=\frac{2c_5(1+\frac{\lambda M^2}{d})^4M^4}{\big(E(\rho^2)\big)^2}$, $c_7(d)=\frac{2c_3(1+\frac{\lambda M^2}{d})^4M^6}{d^2\big(E(\rho^2)\big)^3}$, $c_8(d)=\frac{2c_3(1+\frac{\lambda M^2}{d})^4M^4}{d^2\big(E(\rho^2)\big)^2}$ and
\[
c_9(d)=\frac{1}{d}\frac{d}{\lambda E(\rho^2)}(1+\frac{\lambda M^2}{d})^2=\frac{\big(1+\frac{\lambda M^2}{d}\big)^2}{\lambda E(\rho^2)}.
\]
Since $\lambda>\frac{1}{E(\rho^2)}$, there exists $c_{10}\in (0,1)$  which does not depend on $d$ that
\[
\widehat{c}_5<\frac{3M^4c_5}{\big(E(\rho^2)\big)^2} \text{~and~}c_9(d)\leq c_{10}
\]
for sufficiently large $d$. For sufficiently large $d$,
\[
c_8(d)\leq \frac{1}{2} \text{~and~}\frac{c_7(d)c_8(d)}{1-c_8(d)}\frac{c_{10}}{1-c_{10}}\leq \frac{1}{10}
\]
since $\lim_{d\rightarrow+\infty}c_7(d)=\lim_{d\rightarrow+\infty}c_8(d)=0$. As a result, by Equation \eqref{equ 4.4.13},
\begin{align*}
&\mathbb{E}_d\Big(R(x,y)\Big|\tau_{x,y}<+\infty\Big)\\
&\leq \frac{1}{1-c_8(d)}\frac{3M^4c_5}{\big(E(\rho^2)\big)^2}\sum_{T=0}^{+\infty}\big(c_7(d)\frac{c_8(d)}{1-c_8(d)}\frac{c_{10}}{1-c_{10}}\big)^T\\
&\leq \frac{6M^4c_5}{\big(E(\rho^2)\big)^2}\sum_{T=0}^{+\infty}(\frac{1}{10})^T=\frac{20M^4c_5}{3\big(E(\rho^2)\big)^2}
\end{align*}
for sufficiently large $d$. Let
\[
c_2=\frac{20M^4c_5}{3\big(E(\rho^2)\big)^2}
\]
and the proof is complete.

\qed

At the end of this subsection, we give the proof of Lemma \ref{lemma 4.4.3}.

\proof[Proof of Lemma \ref{lemma 4.4.3}]

For each $m\geq 1$ and each $x\in \mathbb{Z}^d_+$, we define
\[
L_m(x)=\big\{\vec{x}=(x_0,x_1,\ldots,x_m):x_0=x,x_i\rightarrow x_{i+1}\text{~for all~}0\leq i\leq m-1\big\}
\]
as the set of oriented paths starting at $x$ with length $m$.

For each $\vec{x}=(x,x_1,\ldots,x_m)\in L_m(x)$, we denote by $\pi_{\vec{x}}$ the event that $\widetilde{U}(x_i,x_{i+1})<\widetilde{Y}(x_i)$ for all $0\leq i\leq m-1$, then for each $x\in \mathbb{Z}^+_d$,
\begin{align}\label{equ 4.4.14}
P_{\lambda,d}\Big(\pi_{\vec{x}}\Big)&=E_{\mu_d}\Big(\prod_{i=0}^{m-1}P_{\lambda,\omega}\big(\widetilde{U}(x_i,x_{i+1})<\widetilde{Y}(x_i)\big)\Big) \notag\\
&=E_{\mu_d}\Big(\prod_{i=0}^{m-1}\frac{\frac{\lambda}{d}\rho(x_i)\rho(x_{i+1})}{1+\frac{\lambda}{d}\rho(x_i)\rho(x_{i+1})}\Big)\\
&=E\Big(\prod_{i=0}^{m-1}\frac{\frac{\lambda}{d}\rho_i\rho_{i+1}}{1+\frac{\lambda}{d}\rho_i\rho_{i+1}}\Big),
\notag
\end{align}
where $\rho_0,\ldots,\rho_m$ are independent copies of $\rho$. For
$x,y$ that $\|x\|\leq \|y\|$, $m\geq \|y\|-\|x\|$ and
$\vec{x}=(x,x_1,\ldots,x_m)\in L_m(x), \vec{y}=(y,y_1,\ldots,y_m)\in
L_m(y)$,
\begin{align}\label{equ 4.4.15}
&P_{\lambda,d}\Big(\pi_{\vec{x}}\bigcap \pi_{\vec{y}}\Big)=E_{\lambda,d}\Bigg(\prod_{i=0}^{\|y\|-\|x\|}P_{\lambda,\omega}\big(\widetilde{U}(x_i,x_{i+1})<\widetilde{Y}(x_i)\big) \\
&\times\prod_{l=1}^{m-\|y\|+\|x\|} G(x_{l-\|y\|+\|x\|-1},y_{l-1};x_{l-\|y\|+\|x\|},y_l)\notag\\
&\times \prod_{j=m-\|y\|+\|x\|}^{m-1}P_{\lambda,\omega}\big(\widetilde{U}(y_i,y_{j+1})<\widetilde{Y}(y_{j+1})\big)\Bigg), \notag
\end{align}
where
\[
G(x,y;u,v)=P_{\lambda,\omega}\big(\widetilde{U}(x,u)<\widetilde{Y}(x), \widetilde{U}(y,v)<\widetilde{Y}(y)\big)
\]
for $x\rightarrow u$ and $y\rightarrow v$.

By direct calculation, for $x,y$ that $\|x\|=\|y\|$,
\begin{equation}\label{equ 4.4.16}
G(x,y;u,v)
\begin{cases}
=\frac{\frac{\lambda}{d}\rho(x)\rho(u)}{1+\frac{\lambda}{d}\rho(x)\rho(u)} &\text{~if~} x=y \text{~and~}u=v,\\
\leq \frac{2\frac{\lambda^2}{d^2}\rho^2(x)\rho(u)\rho(v)}{[1+\frac{\lambda}{d}\rho(x)\rho(u)][1+\frac{\lambda}{d}\rho(x)\rho(v)]}
&\text{~if~} x=y \text{~and~}u\neq v,\\
=\frac{\frac{\lambda^2}{d^2}\rho(x)\rho(y)\rho^2(u)}{[1+\frac{\lambda}{d}\rho(x)\rho(u)][1+\frac{\lambda}{d}\rho(y)\rho(u)]}
&\text{~if~} x\neq y \text{~and~}u=v\\
=\frac{\frac{\lambda^2}{d^2}\rho(x)\rho(y)\rho(u)\rho(v)}{[1+\frac{\lambda}{d}\rho(x)\rho(u)][1+\frac{\lambda}{d}\rho(y)\rho(v)]}
&\text{~if~} x\neq y \text{~and~}u\neq v.
\end{cases}
\end{equation}
According to the definition of the SIR model, for given $A\subseteq
\mathbb{Z}^d_+$, $m\geq 1$ and $x\in A$, if $\pi_{\vec{x}}$ occurs
for some $\vec{x}=(x,x_1,\ldots,x_m)\in L_m(x)$, then
\[
x_m\in \bigcup_{t\geq 0}I_t^A.
\]
As a result, on the event $\bigcap_{m=1}^{+\infty}\bigcup_{x\in
A}\bigcup_{\vec{x}\in L_m(x)}\pi_{\vec{x}}$, there are infinite many
vertices have ever been infected and hence
\begin{equation}\label{equ 4.4.17}
P_{\lambda,d}\big(I_t^A\neq \emptyset,\forall~t\geq 0\big)\geq
P_{\lambda,d}\Big(\bigcap_{m=1}^{+\infty}\bigcup_{x\in
A}\bigcup_{\vec{x}\in L_m(x)}\pi_{\vec{x}}\Big).
\end{equation}
We use $\chi_{\vec{x}}$ to denote the indicator function of
$\pi_{\vec{x}}$, then by the Cauchy-Schwartz's inequality and the
dominated convergence theorem,
\begin{align*}
&P_{\lambda,d}\Big(\bigcap_{m=1}^{+\infty}\bigcup_{x\in
A}\bigcup_{\vec{x}\in L_m(x)}\pi_{\vec{x}}\Big)\\
&\geq \lim_{m\geq
1}P_{\lambda,d}\big(\bigcup_{x\in A}\bigcup_{\vec{x}\in
L_m(x)}\pi_{\vec{x}}\big)\\
&=\lim_{m\geq 1}P_{\lambda,d}\big(\sum_{x\in A}\sum_{\vec{x}\in
L_m(x)}\chi_{\vec{x}}>0\big)\\
&\geq \limsup_{m\rightarrow+\infty}
\frac{\Bigg[E_{\lambda,d}\Big(\sum_{x\in A}\sum_{\vec{x}\in
L_m(x)}\chi_{\vec{x}}\Big)\Bigg]^2}{E_{\lambda,d}\Bigg[\Big(\sum_{x\in
A}\sum_{\vec{x}\in L_m(x)}\chi_{\vec{x}}\Big)^2\Bigg]}\\
&=\limsup_{m\rightarrow+\infty}\frac{\Bigg[\sum\limits_{x\in
A}\sum\limits_{\vec{x}\in
L_m(x)}P_{\lambda,d}\big(\pi_{\vec{x}}\big)\Bigg]^2}{\sum\limits_{x\in
A}\sum\limits_{y\in A}\sum\limits_{\vec{x}\in
L_m(x)}\sum\limits_{\vec{y}\in
L_m(y)}P_{\lambda,d}\big(\pi_{\vec{x}}\bigcap\pi_{\vec{y}}\big)}.
\end{align*}
By Equation \eqref{equ 4.4.14}, for given $m\geq 1$ and $\vec{x}\in
L_m(x)$, $P_{\lambda,d}\big(\pi_{\vec{x}}\big)$ does not depend on
the choice of $\vec{x}$ and $x$. Therefore, according to the fact
that $|L_m(x)|=d^m$,
\begin{align}\label{equ 4.4.18}
&P_{\lambda,d}\Big(\bigcap_{m=1}^{+\infty}\bigcup_{x\in
A}\bigcup_{\vec{x}\in L_m(x)}\pi_{\vec{x}}\Big)\\
&\geq
\frac{1}{\liminf\limits_{m\rightarrow+\infty}\frac{1}{|A|^2}\sum\limits_{x\in
A}\sum\limits_{y\in A}\sum\limits_{\vec{x}\in
L_m(x)}\sum\limits_{\vec{y}\in
L_m(y)}\frac{1}{d^{2m}}\frac{P_{\lambda,d}\big(\pi_{\vec{x}}\bigcap\pi_{\vec{y}}\big)}
{P_{\lambda,d}\big(\pi_{\vec{x}}\big)P_{\lambda,d}\big(\pi_{\vec{y}}\big)}}.
\notag
\end{align}
We use $\vec{\vartheta}_m^x$ to denote the random path
$(\vartheta_0^x,\ldots,\vartheta_m^x)$ while use $\vec{\nu}_m^y$ to
denote the path $(\nu_0^y,\ldots,\nu_m^y)$, then by Equation
\eqref{equ 4.4.18},
\begin{equation}\label{equ 4.4.19}
P_{\lambda,d}\Big(\bigcap_{m=1}^{+\infty}\bigcup_{x\in
A}\bigcup_{\vec{x}\in L_m(x)}\pi_{\vec{x}}\Big)\geq
\frac{1}{\liminf\limits_{m\rightarrow+\infty}\frac{1}{|A|^2}\sum\limits_{x\in
A}\sum\limits_{y\in
A}\mathbb{E}_d\Big(\frac{P_{\lambda,d}(\pi_{_{\vec{\vartheta}_m^x}}\bigcap\pi_{_{\vec{\nu}_m^y}})}
{P_{\lambda,d}(\pi_{_{\vec{\vartheta}_m^x}})P_{\lambda,d}(\pi_{_{\vec{\nu}_m^y}})}\Big)}.
\end{equation}

For $m$ sufficiently large and $x,y$ that $\|x\|\leq \|y\|$, we
bound
$\frac{P_{\lambda,d}(\pi_{_{\vec{\vartheta}_m^x}}\bigcap\pi_{_{\vec{\nu}_m^y}})}
{P_{\lambda,d}(\pi_{_{\vec{\vartheta}_m^x}})P_{\lambda,d}(\pi_{_{\vec{\nu}_m^y}})}$
from above according to the following procedure. For the denominator
\[
P_{\lambda,d}(\pi_{_{\vec{\vartheta}_m^x}})P_{\lambda,d}(\pi_{_{\vec{\nu}_m^y}})=
E_{\mu_d}\Big(\prod_{i=0}^{m-1}\frac{\frac{\lambda}{d}\rho(\vartheta_i)\rho(\vartheta_{i+1})}{1+\frac{\lambda}{d}\rho(\vartheta_i)\rho(\vartheta_{i+1})}\Big)
E_{\mu_d}\Big(\prod_{i=0}^{m-1}\frac{\frac{\lambda}{d}\rho(\nu_i)\rho(\nu_{i+1})}{1+\frac{\lambda}{d}\rho(\nu_i)\rho(\nu_{i+1})}\Big),
\]
if $l\geq \|y\|-\|x\|$ satisfies that $\vartheta_l^x=\nu_{l-\|y\|+\|x\|}^y$, then
\[
\begin{cases}
&\frac{\frac{\lambda}{d}\rho(\vartheta_l^x)\rho(\vartheta_{l+1}^x)}{1+\frac{\lambda}{d}\rho(\vartheta_l^x)\rho(\vartheta_{l+1}^x)}\geq \frac{\frac{\lambda}{d}\rho(\vartheta_l^x)\rho(\vartheta_{l+1}^x)}{1+\frac{\lambda}{d}M^2}, \\ &\frac{\frac{\lambda}{d}\rho(\nu_{l-\|y\|+\|x\|}^y)\rho(\nu_{l-\|y\|+\|x\|+1}^y)}{1+\frac{\lambda}{d}\rho(\nu_{l-\|y\|+\|x\|}^y)\rho(\nu_{l-\|y\|+\|x\|+1}^y)}\geq \frac{\frac{\lambda}{d}\rho(\nu_{l-\|y\|+\|x\|}^y)\rho(\nu_{l-\|y\|+\|x\|+1}^y)}{1+\frac{\lambda}{d}M^2},\\
&\frac{\frac{\lambda}{d}\rho(\vartheta_{l-1}^x)\rho(\vartheta_{l}^x)}{1+\frac{\lambda}{d}\rho(\vartheta_{l-1}^x)\rho(\vartheta_{l}^x)}\geq \frac{\frac{\lambda}{d}\rho(\vartheta_{l-1}^x)\rho(\vartheta_{l}^x)}{1+\frac{\lambda}{d}M^2},\\
&\frac{\frac{\lambda}{d}\rho(\nu_{l-\|y\|+\|x\|-1}^y)\rho(\nu_{l-\|y\|+\|x\|}^y)}{1+\frac{\lambda}{d}\rho(\nu_{l-\|y\|+\|x\|-1}^y)\rho(\nu_{l-\|y\|+\|x\|}^y)}\geq \frac{\frac{\lambda}{d}\rho(\nu_{l-\|y\|+\|x\|-1}^y)\rho(\nu_{l-\|y\|+\|x\|})}{1+\frac{\lambda}{d}M^2}.
\end{cases}
\]
For the numerator $P_{\lambda,d}(\pi_{_{\vec{\vartheta}_m^x}}\bigcap\pi_{_{\vec{\nu}_m^y}})$ with expression given by Equation \eqref{equ 4.4.15},
\[
\begin{cases}
& G(\vartheta_l^x,\nu_{l-\|y\|+\|x\|}^y;\vartheta_{l+1}^x,\nu_{l-\|y\|+\|x\|+1}^y)\leq \frac{\lambda}{d}\rho(\vartheta_l^x)\rho(\vartheta_{l+1}^x) \\
& \text{~if~}\vartheta_l^x=\nu_{l-\|y\|+\|x\|}^y\text{~and~}\vartheta_{l+1}^x=\nu_{l-\|y\|+\|x\|+1}^y,\\
& G(\vartheta_l^x,\nu_{l-\|y\|+\|x\|}^y;\vartheta_{l+1}^x,\nu_{l-\|y\|+\|x\|+1}^y)\leq \frac{2\lambda^2}{d^2}\rho^2(\vartheta_l^x)\rho(\vartheta_{l+1}^x)\rho(\nu_{l-\|y\|+\|x\|+1}^y) \\
& \text{~if~}\vartheta_l^x=\nu_{l-\|y\|+\|x\|}^y \text{~and~} \vartheta_{l+1}^x\neq\nu_{l-\|y\|+\|x\|+1}^y,\\
& G(\vartheta_l^x,\nu_{l-\|y\|+\|x\|}^y;\vartheta_{l+1}^x,\nu_{l-\|y\|+\|x\|+1}^y)\leq \frac{\lambda^2}{d^2}\rho(\vartheta_l^x)\rho(\nu_{l-\|y\|+\|x\|}^y)\rho^2(\vartheta_{l+1}^x) \\
& \text{~if~}\vartheta_l^x\neq \nu_{l-\|y\|+\|x\|}^y \text{~and~} \vartheta_{l+1}^x=\nu_{l-\|y\|+\|x\|+1}^y.
\end{cases}
\]
According to the aforesaid inequalities, $\frac{P_{\lambda,d}(\pi_{_{\vec{\vartheta}_m^x}}\bigcap\pi_{_{\vec{\nu}_m^y}})}
{P_{\lambda,d}(\pi_{_{\vec{\vartheta}_m^x}})P_{\lambda,d}(\pi_{_{\vec{\nu}_m^y}})}$ is bounded from above by an upper bound $R_m(x,y)$. According to our assumption of the independence between the exponential times, the expression of $R_m(x,y)$ can be simplified by canceling common factors in the numerator and denominator. For example, if $l<k$ that $\vartheta_l^x=\nu_{l-\|y\|+\|x\|}^y$ and $\vartheta_{k}^x=\nu_{k-\|y\|+\|x\|}^y$ while $\vartheta_j^x\neq \nu_{j-\|y\|+\|x\|}^y$ for any $l<j<k$, then both the numerator and denominator have the factor
\[
\Bigg(E\big(\frac{\prod_{i=1}^{k-l-1}\rho_i^2}{\prod_{i=1}^{k-l-2}(1+\rho_i\rho_{i+1})}\big)\Bigg)^2
\]
that can be canceled, where $\rho_1,\ldots,\rho_{k-1-l}$ are independent copies of $\rho$. As a result, it is not difficult to check that
\[
\lim_{m\rightarrow+\infty}R_m(x,y)=R(x,y)
\]
and hence
\begin{align}\label{equ 4.4.20}
P_{\lambda,d}\Big(\bigcap_{m=1}^{+\infty}\bigcup_{x\in
A}\bigcup_{\vec{x}\in L_m(x)}\pi_{\vec{x}}\Big)&\geq
\frac{1}{\lim\limits_{m\rightarrow+\infty}\frac{1}{|A|^2}\sum\limits_{x\in
A}\sum\limits_{y\in
A}\mathbb{E}_d\big(R_m(x,y)\big)} \notag\\
&=\frac{1}{\frac{1}{|A|^2}\sum\limits_{x\in
A}\sum\limits_{y\in
A}\mathbb{E}_d\big(R(x,y)\big)}
\end{align}
according to Equation \eqref{equ 4.4.19}. Lemma \ref{lemma 4.4.3} follows directly from Equations \eqref{equ 4.4.17} and \eqref{equ 4.4.20}. 

\qed

\section{Proof of Equation \eqref{equ 2.4}}\label{section five}
In this section we give the proof of Equation \eqref{equ 2.4}. We
still assume that the vertex weight $\rho$ satisfies \eqref{equ
4.1}. The assumption is without loss of generality according to the
following analysis. For general $\rho$ not satisfying \eqref{equ
4.1}, we let
\[
\widehat{\rho}_m
=
\begin{cases}
\rho &\text{~if~}\rho\geq \frac{1}{m},\\
\frac{1}{m} &\text{~if~}\rho<\frac{1}{m},
\end{cases}
\]
then $\widehat{\rho}_m\geq \rho$ and
$\lim_{m\rightarrow+\infty}\widehat{\rho}_m=\rho$. Therefore,
\[
P_{\lambda,d,\rho}\big(C_t^O\neq \emptyset,\forall~t\geq 0\big)\leq
P_{\lambda,d,\widehat{\rho}_m}\big(C_t^O\neq \emptyset,\forall~t\geq
0\big).
\]
If Equation \eqref{equ 2.4} holds under assumption \eqref{equ 4.1},
which $\widehat{\rho}_{m}$ satisfies, then
\[
\limsup_{d\rightarrow+\infty}P_{\lambda,d,\rho}\big(C_t^O\neq
\emptyset,\forall~t\geq 0\big)\leq
\limsup_{d\rightarrow+\infty}P_{\lambda,d,\widehat{\rho}_m}\big(C_t^O\neq
\emptyset,\forall~t\geq 0\big)\leq
E\big(\frac{\lambda\widehat{\rho}_m\widehat{\theta}_m}{1+\lambda\widehat{\rho}_m\widehat{\theta}_m}\big),
\]
where $\widehat{\theta}_m$ satisfies
\[
E\big(\frac{\lambda\widehat{\rho}_m^2}{1+\lambda\widehat{\rho}_m\widehat{\theta}_m}\big)=1
\]
and it is easy to check that
$\lim_{m\rightarrow+\infty}\widehat{\theta}_m=\theta$. Let
$m\rightarrow+\infty$, then Equation \eqref{equ 2.4} holds for
general $\rho$.

For each $n\geq 0$, we define
\[
\beta_n=\big\{x\in \mathbb{Z}_+^d:\|x\|=n\text{~and~}x\in
\bigcup_{t\geq 0}C_t^O\big\}
\]
as the vertices that with $l_1$ norm $n$ and have ever been infected
with respect to the contact process with $O$ as the unique initially
infected vertex.

Since each infected vertex waits for an exponential with rate one to
become healthy, the infected vertices never die out when and only
when there are infinite many vertices that have ever been infected.
Furthermore, since $x$ infects $y$ only if $x\rightarrow y$,
\begin{equation}\label{equ 5.01}
\big\{C_t^O\neq \emptyset,\forall~t\geq 0\big\}=\big\{\beta_n\neq
\emptyset \text{~for all~}n\geq 0\big\}.
\end{equation}

The proof of Equation \eqref{equ 2.4} relies heavily on  Equation
\eqref{equ 5.01} and the following two lemmas.

\begin{lemma}\label{lemma 5.1}
Let $\{W_n\}_{n\geq 0}$ be the branching process with random vertex
weights defined as in Section \ref{section three} and $\sigma\in
(0,\frac{1}{10\log(\lambda M^2)})$, then
\[
\liminf_{d\rightarrow+\infty}\widehat{P}_{\lambda,d}\big(W_{\lfloor\sigma\log
d\rfloor}=\emptyset\big)\geq
E\big(\frac{1}{1+\lambda\rho\theta}\big),
\]
where $\widehat{P}_{\lambda,d}$ is the annealed measure of the
branching process defined as in Section \ref{section three}.

\end{lemma}

\begin{lemma}\label{lemma 5.2}
Let $\{V_n\}_{n\geq 0}$ be defined as in Section \ref{section four}
and $\sigma\in (0,\frac{1}{10\log(\lambda M^2)})$, then
\[
\lim_{d\rightarrow+\infty}\Big[P_{\lambda,d}\big(\beta_{\lfloor\sigma\log
d\rfloor}=\emptyset\big)-P_{\lambda,d}\big(V_{\lfloor\sigma\log
d\rfloor}=\emptyset\big)\Big]=0.
\]
\end{lemma}

The proof of Lemma \ref{lemma 5.1} is given in Subsection
\ref{subsection 5.1}. The core idea of the proof is to show that the
branching process survives with high probability conditioned on
$W_{\lfloor\sigma\log d\rfloor}\neq \emptyset$. The proof of Lemma
\ref{lemma 5.2} is given in Subsection \ref{subsection 5.2}. The
core idea of the proof is to construct a coupling of
$\{\beta_n\}_{n\geq 0}$ and $\{V_n\}_{n\geq 0}$ such that
$\beta_{\lfloor\sigma\log d\rfloor}=V_{\lfloor\sigma\log d\rfloor}$
with high probability. Now we show how to utilize Lemmas \ref{lemma
5.1} and \ref{lemma 5.2} to prove Equation \eqref{equ 2.4}.

\proof[Proof of Equation \eqref{equ 2.4}]

We couple $\{W_n\}_{n\geq 0}$ and $\{V_n\}_{n\geq 0}$ under the same
probability space as what we have done in Subsection \ref{subsection
4.1}. Recalling that we define $B(d)$ as the event that the coupling
is successful at step $m$ for all $m\leq \lfloor\sigma\log
d\rfloor$, then
\[
V_{\lfloor\sigma\log d\rfloor}=W_{\lfloor\sigma\log d\rfloor}
\]
on the event $B(d)$. Therefore, by Lemma \ref{lemma 4.1.2},
\[
\Big|\widehat{P}_{\lambda,d}\big(W_{\lfloor\sigma\log
d\rfloor}=\emptyset\big)-P_{\lambda,d}\big(V_{\lfloor\sigma\log
d\rfloor}=\emptyset\big)\Big|\leq
2P_{\lambda,d}\big(B(d)^c\big)\rightarrow 0
\]
as $d\rightarrow+\infty$ and hence
\[
\liminf_{d\rightarrow+\infty}P_{\lambda,d}\big(V_{\lfloor\sigma\log
d\rfloor}=\emptyset\big)\geq
E\big(\frac{1}{1+\lambda\rho\theta}\big)
\]
according to Lemma \ref{lemma 5.1}. Then, by Lemma \ref{lemma 5.2},
\[
\liminf_{d\rightarrow+\infty}P_{\lambda,d}\big(\beta_{\lfloor\sigma\log
d\rfloor}=\emptyset\big)\geq
E\big(\frac{1}{1+\lambda\rho\theta}\big)
\]
and hence
\begin{equation}\label{equ 5.02}
\liminf_{d\rightarrow+\infty}P_{\lambda,d}\big(\beta_n=\emptyset\text{~for
some~}n\geq 0\big)\geq E\big(\frac{1}{1+\lambda\rho\theta}\big).
\end{equation}
By Equation \eqref{equ 5.02},
\begin{equation}\label{equ 5.03}
\limsup_{d\rightarrow+\infty}P_{\lambda,d}\big(\beta_n\neq\emptyset\text{~for
all~}n\geq 0\big)\leq
E\big(\frac{\lambda\rho\theta}{1+\lambda\rho\theta}\big).
\end{equation}

Equation \eqref{equ 2.4} follows from Equations \eqref{equ 5.01} and
\eqref{equ 5.03} directly.

\qed

\subsection{Proof of Lemma \ref{lemma 5.1}}\label{subsection 5.1}
In this subsection, we give the proof of Lemma \ref{lemma 5.1}.
First we introduce some notations and definitions. For sufficiently
large $d$, let $N(d)=\log(\log d)$ defined as in Section
\ref{section four}. For each $x\in \mathbb{T}^d$, we give the $d$
sons of $x$ an order $x(1),x(2),\ldots,x(d)$. Then, we define

(1) $\widehat{W}_0=\Upsilon$.

(2) For each $n\geq 0$,
\begin{align*}
\widehat{W}_{n+1}=&\big\{y:\text{there exists~}x\in
\widehat{W}_{n}\text{~that~}y=x(i)\\
&\text{~for some~}i\leq d-\lfloor\frac{d}{N(d)}\rfloor
\text{~and~}U(x,y)<Y(x)\big\}.
\end{align*}

According to the aforesaid definition, $\{\widehat{W}_n\}_{n\geq 0}$
is a branching process with random vertex weights on a subtree of
$\mathbb{T}^d$ while this subtree is isomorphic to
$\mathbb{T}^{d-\lfloor\frac{d}{N(d)}\rfloor}$. For each $n\geq 0$,
$\widehat{W}_n\subseteq W_n$. The following lemma is crucial for us
to prove Lemma \ref{lemma 5.1}.

\begin{lemma}\label{lemma 5.1.1}
For $\sigma\in (0,\frac{1}{10\log (\lambda M^2)})$,
\[
\lim_{d\rightarrow+\infty}\widehat{P}_{\lambda,d}\Big(W_n\neq
\emptyset\text{~for all~}n\geq 0\Big|\widehat{W}_{\lfloor\sigma\log
d\rfloor}\neq \emptyset\Big)=1.
\]
\end{lemma}

We give the proof of Lemma \ref{lemma 5.1.1} at the end of this
subsection. Now we show how to utilize Lemma \ref{lemma 5.1.1} to
prove Lemma \ref{lemma 5.1}.

\proof[Proof of Lemma \ref{lemma 5.1}]

By conditional probability formula,
\begin{align*}
&\widehat{P}_{\lambda,d}\big(W_n\neq \emptyset \text{~for all~}n\geq
0\big) \\
&\geq \widehat{P}_{\lambda,d}\Big((W_n\neq \emptyset\text{~for
all~}n\geq 0\Big|\widehat{W}_{\lfloor\sigma\log d\rfloor}\neq
\emptyset\Big)\widehat{P}_{\lambda,d}\big(\widehat{W}_{\lfloor\sigma\log
d\rfloor}\neq \emptyset\big).
\end{align*}
Then, by Lemmas \ref{lemma 3.1} and \ref{lemma 5.1.1},
\[
\limsup_{d\rightarrow+\infty}
\widehat{P}_{\lambda,d}\big(\widehat{W}_{\lfloor\sigma\log
d\rfloor}\neq \emptyset\big)\leq
E\big(\frac{\lambda\rho\theta}{1+\lambda\rho\theta}\big)
\]
and hence
\begin{equation}\label{equ 5.1.1}
\liminf_{d\rightarrow+\infty}
\widehat{P}_{\lambda,d}\big(\widehat{W}_{\lfloor\sigma\log
d\rfloor}=\emptyset\big)\geq
E\big(\frac{1}{1+\lambda\rho\theta}\big)
\end{equation}
for any $\sigma\in (0,\frac{1}{10\log(\lambda M^2)})$.

For given $\lambda>\frac{1}{E(\rho^2)}$ and $\sigma\in
(0,\frac{1}{10\log (\lambda M^2)})$, we choose arbitrary
$\widehat{\lambda}\in (\lambda,+\infty)$ and $\widehat{\sigma}\in
(0,\sigma)$. For sufficiently large $d$, we define
\[
\widehat{d}=\inf\{k:k-\lfloor\frac{k}{N(k)}\rfloor\geq d\},
\]
then it is easy to check that
$\lim_{d\rightarrow+\infty}\frac{\widehat{d}}{d}=1$ and hence
\begin{equation}\label{equ 5.1.2}
\frac{\widehat{\lambda}}{\widehat{d}}\geq
\frac{\lambda}{d}\text{~while~}\widehat{\sigma}\log \widehat{d}\leq
\sigma\log d
\end{equation}
for sufficiently large $d$. As we have introduced,
$\{\widehat{W}_n\}_{n\geq 0}$ on $\mathbb{T}^{\widehat{d}}$ can be
identified with $\{W_n\}_{n\geq 0}$ on
$\mathbb{T}^{\widehat{d}-\lfloor\frac{\widehat{d}}{N(\widehat{d})}\rfloor}$
with a scaling of the infection rate $\lambda$. As a result, by
Equation \eqref{equ 5.1.2},
\begin{equation}\label{equ 5.1.3}
\widehat{P}_{\widehat{\lambda},\widehat{d}}\big(\widehat{W}_{\lfloor\widehat{\sigma}\log
\widehat{d}\rfloor}=\emptyset\big)\leq
\widehat{P}_{\lambda,d}\big(W_{\lfloor\sigma\log
d\rfloor}=\emptyset\big)
\end{equation}
for sufficiently large $d$. By Equations \eqref{equ 5.1.1} and
\eqref{equ 5.1.3},
\begin{equation}\label{equ 5.1.4}
\liminf_{d\rightarrow+\infty}\widehat{P}_{\lambda,d}\big(W_{\lfloor\sigma\log
d\rfloor}=\emptyset\big)\geq \liminf_{d\rightarrow+\infty}
\widehat{P}_{\widehat{\lambda},\widehat{d}}\big(\widehat{W}_{\lfloor\widehat{\sigma}\log
\widehat{d}\rfloor}=\emptyset\big)\geq
E\big(\frac{1}{1+\widehat{\lambda}\rho\widehat{\theta}}\big),
\end{equation}
where $\widehat{\theta}$ satisfies
\[
E\big(\frac{\widehat{\lambda}\rho^2}{1+\widehat{\lambda}\rho\widehat{\theta}}\big)=1
\]
and it is easy to check that
$\lim_{\widehat{\lambda}\rightarrow\lambda}\widehat{\theta}=\theta$.

Let $\widehat{\lambda}\rightarrow\lambda$, then Lemma \ref{lemma
5.1} follows directly from Equation \eqref{equ 5.1.4}.

\qed

At the end of this subsection we give the proof of Lemma \ref{lemma
5.1.1}.

\proof[Proof of Lemma \ref{lemma 5.1.1}]

We define
\begin{align*}
\widehat{D}=\big\{&y:\text{~there exists~}x\in
\bigcup_{m=0}^{\lfloor\sigma\log
d\rfloor}\widehat{W}_m\\
&\text{~that~}y=x(i)\text{~for some~}i\geq
d-\lfloor\frac{d}{N(d)}\rfloor+1\text{~and~}U(x,y)<Y(x)\big\}.
\end{align*}
It is easy to check that $u\in \bigcup_{m\geq 0}W_m$ and $\rho(u)>0$
for any $u\in \widehat{D}$. According to nearly the same analysis as
that in the proof of Lemma \ref{lemma 4.2},
\begin{equation}\label{equ 5.1.5}
\lim_{d\rightarrow+\infty}\widehat{P}_{\lambda,d}\Big(|\widehat{D}|\geq
K(d)\Big|\widehat{W}_{\lfloor\sigma\log d\rfloor}\neq
\emptyset\Big)=1,
\end{equation}
where $K(d)=\lfloor\frac{\sqrt{\log d}}{N(d)}\rfloor$ defined as in
Section \ref{section four}.

For any $u\in \widehat{D}$, we denote by $\mathbb{T}_u$ the subtree
of $\mathbb{T}^d$ rooted at $u$ and consisted of $u$ and its
descendants. For $u\in D$, $\rho(u)>0$ and hence
$\rho(u)\geq\epsilon$ by assumption \eqref{equ 4.1}. Then, the SIR
model confined on $\mathbb{T}_u$ with $u$ as the unique initially
infected vertex has survival probability at least
\begin{equation}\label{equ 5.1.6}
1-F_d(\epsilon),
\end{equation}
where $F_d(\epsilon)$ is defined as in Section \ref{section three}.
It is easy to check that $\mathbb{T}_u\bigcap
\mathbb{T}_v=\emptyset$ for any $u,v\in D,u\neq v$, then by Equation
\eqref{equ 5.1.6},
\begin{align}\label{equ 5.1.7}
&P_{\lambda,d}\Big(W_n\neq \emptyset\text{~for all~}n\geq
0\Big||\widehat{D}|\geq K(d)\Big)\\
&\geq
1-\Big[1-(1-F_d(\epsilon))\Big]^{K(d)}=1-F_d(\epsilon)^{K(d)},\notag
\end{align}
since all the vertices in $\widehat{D}$ have ever been infected.

By Lemmas \ref{lemma 3.1} and \ref{lemma 3.2},
\[
\lim_{d\rightarrow+\infty}F_d(\epsilon)=\frac{1}{1+\lambda
\epsilon\theta}
\]
and hence $\lim_{d\rightarrow+\infty}1-F_d(\epsilon)^{K(d)}=1$. Then
by Equation \eqref{equ 5.1.7},
\begin{equation}\label{equ 5.1.8}
\lim_{d\rightarrow+\infty}P_{\lambda,d}\Big(W_n\neq
\emptyset\text{~for all~}n\geq 0\Big||\widehat{D}|\geq K(d)\Big)=1.
\end{equation}
Lemma \ref{lemma 5.1.1} follows from Equations \eqref{equ 5.1.5} and
\eqref{equ 5.1.8} directly.

\qed

\subsection{Proof of Lemma \ref{lemma 5.2}}\label{subsection 5.2}
In this subsection we give the proof of Lemma \ref{lemma 5.2}. First we couple $\{\beta_n\}_{n\geq 0}$ and $\{V_n\}_{n\geq 0}$ under the same probability space.
Let $\{\widetilde{Y}(x)\}_{x\in \mathbb{Z}^d_+}$ and $\{\widetilde{U}(x,y)\}_{x\in \mathbb{Z}^d_+,x\rightarrow y}$ be defined as in Section \ref{section four}, then $\{V_n\}_{n\geq 0}$ is defined as in Section \ref{section four} according to the values of $\{\widetilde{Y}(x)\}_{x\in \mathbb{Z}^d_+}$ and $\{\widetilde{U}(x,y)\}_{x\in \mathbb{Z}^d_+,x\rightarrow y}$. For any $x,y\in \mathbb{Z}^d_+,x\rightarrow y$, let $\widetilde{U}_2(x,y)$ be an independent copy of $\widetilde{U}(x,y)$ under the quenched measure. We assume that all these exponential times are independent under the quenched measure. For the contact process, we let $\widetilde{Y}(x)$ be the time $x$ waits for to become healthy after the first moment when $x$ is infected. We let $\widetilde{U}(x,y)$ be the time $x$ waits for to infect $y$ after the first moment when $x$ is infected. If $\widetilde{U}(x,y)<\widetilde{Y}(x)$, then after the first infection from $x$ to $y$, $x$ waits for $\widetilde{U}_2(x,y)$ units of time to infect $y$ again, i.e., $x$ infects $y$ at least twice before becoming healthy when $\widetilde{U}(x,y)+\widetilde{U}_2(x,y)<\widetilde{Y}(x)$. Following the above definitions, $\{V_n\}_{n\geq 0}$ and $\{\beta_n\}_{n\geq 0}$ are coupled under the same probability space and it is obviously that $V_n\subseteq \beta_n$ for each $n\geq 0$.

Let $J(d)$ be defined as in Section \ref{section four}, i.e., the event that $\widetilde{U}(x,y)>\widetilde{Y}(x)$ and $\widetilde{U}(z,y)>\widetilde{Y}(z)$ for any $x,y,z$ that $x,z\in \bigcup_{m=0}^{\lfloor\sigma\log d\rfloor}V_m$ and $x,z\rightarrow y$. On the event $J(d)$, if $V_{\lfloor\sigma\log d\rfloor}\neq \beta_{\lfloor\sigma\log d\rfloor}$, then there must exist repeated infection from some $x$ to $y$ that $x\rightarrow y$ for the contact process, i.e.,
\[
\widetilde{U}(x,y)+\widetilde{U}_2(x,y)<\widetilde{Y}(x).
\]
For each $m\geq 1$, $L_m(O)$ is the set of oriented paths on $\mathbb{Z}^d_+$ starting at $O$ with length $m$ defined as in Section \ref{section four}. For each $\vec{l}:O=l_0\rightarrow l_1\rightarrow l_2\rightarrow \ldots\rightarrow l_m$ in $L_m(O)$, we denote by $\widehat{A}_{\vec{l}}$ the event that $\widetilde{U}(l_i,l_{i+1})<\widetilde{Y}(l_i)$ for all $0\leq i\leq m-2$ and
\[
\widetilde{U}(l_m-1,l_m)+\widetilde{U}_2(l_{m-1},l_m)<\widetilde{Y}(l_{m-1}),
\]
then according to the aforesaid analysis,
\begin{equation}\label{equ 5.2.1}
P_{\lambda,d}\big(V_{\lfloor\sigma\log d\rfloor}\neq \beta_{\lfloor\sigma\log d\rfloor},J(d)\big)\leq \sum_{m=0}^{\lfloor\sigma\log d\rfloor}\sum_{\vec{l}\in L_m}P_{\lambda,d}(\widehat{A}_{\vec{l}}).
\end{equation}
Now we give the proof of Lemma \ref{lemma 5.2}.

\proof[Proof of Lemma \ref{lemma 5.2}]

For each $\vec{l}\in L_m$, since $\widetilde{U}(\cdot,\cdot),\widetilde{U}_2(\cdot,\cdot)$ are exponential times with rate at most $\frac{\lambda M^2}{d}$ while $\widetilde{Y}(\cdot)$ is an exponential time with rate one, it is easy to check that
\[
P_{\lambda,d}(\widehat{A}_{\vec{l}})\leq \big(\frac{\lambda M^2}{d}\big)^{m-1}\big(\frac{\lambda M^2}{d}\big)^2=\frac{\lambda^{m+1}M^{2m+2}}{d^{m+1}}.
\]
Since $|L_m|=d^m$ and $\sigma<\frac{1}{10\log (\lambda M^2)}$, by Equation \eqref{equ 5.2.1},
\begin{equation}\label{equ 5.2.2}
P_{\lambda,d}\big(V_{\lfloor\sigma\log d\rfloor}\neq \beta_{\lfloor\sigma\log d\rfloor},J(d)\big)
\leq \sum_{m=0}^{\lfloor\sigma\log d\rfloor} d^m \frac{\lambda^{m+1}M^{2m+2}}{d^{m+1}}\leq \frac{\lambda^2 M^4d^{-0.9}}{\lambda M^2-1}.
\end{equation}
By Equation \eqref{equ 5.2.2},
\begin{align}\label{equ 5.2.3}
&|P_{\lambda,d}\big(\beta_{\lfloor\sigma\log d\rfloor}=\emptyset\big)-P_{\lambda,d}\big(V_{\lfloor\sigma\log d\rfloor}=\emptyset\big)|\\
&\leq 2P_{\lambda,d}\big(\beta_{\lfloor\sigma\log d\rfloor}\neq V_{\lfloor\sigma\log d\rfloor}\big)\notag\\
&\leq 2P_{\lambda,d}\big(\beta_{\lfloor\sigma\log d\rfloor}\neq V_{\lfloor\sigma\log d\rfloor},J(d)\big)+2P_{\lambda,d}\big(J(d)^c\big)\notag\\
&\leq \frac{2\lambda^2 M^4d^{-0.9}}{\lambda M^2-1}+2P_{\lambda,d}\big(J(d)^c\big).\notag
\end{align}
By Equation \eqref{equ 4.1.4},
\begin{equation}\label{equ 5.2.4}
\lim_{d\rightarrow+\infty}P_{\lambda,d}\big(J(d)^c\big)=0.
\end{equation}
Lemma \ref{lemma 5.2} follows directly from Equations \eqref{equ 5.2.3} and \eqref{equ 5.2.4}.

\qed

\quad

\textbf{Acknowledgments.} The author is grateful to the financial
support from the National Natural Science Foundation of China with
grant number 11501542 and the financial support from Beijing
Jiaotong University with grant number KSRC16006536.

{}
\end{document}